\documentclass[12pt]{article}

\title{Conformal Orbifold Theories \\ and Braided Crossed G-Categories}
\author{Michael M\"uger\thanks{Supported by NWO.} \\
Korteweg-de Vries Institute for Mathematics \\ University of Amsterdam, Netherlands \\ 
email: {\tt mmueger@science.uva.nl}}

\newlength{\dinwidth}
\newlength{\dinmargin}
\usepackage{latexsym, amssymb, amsmath, theorem,tan-gle}
\usepackage{tan-gle}
\input diagrams         

\def\1#1{{\bf #1}}
\def\2#1{{\cal #1}}
\def\3#1{{\sl #1}}
\def\4#1{{\tt #1}}
\def\5#1{{\sf #1}}
\def\6#1{{\mathfrak #1}}
\def\7#1{{\mathbb #1}}

\newcommand{\ba}{\begin{array}}
\newcommand{\ea}{\end{array}}
\newcommand{\bea}{\begin{eqnarray}}
\newcommand{\eea}{\end{eqnarray}}
\newcommand{\bean}{\begin{eqnarray*}}
\newcommand{\eean}{\end{eqnarray*}}

\newcommand{\ve}{\varepsilon}

\newcommand{\impl}{\Rightarrow}
\newcommand{\rarr}{\rightarrow}
\newcommand{\restr}{\upharpoonright}
\newcommand{\ol}{\overline}

\newcommand{\del}{\partial}
\newcommand{\id}{\mathrm{id}}
\newcommand{\mcirc}{\,\circ\,}

\newcommand{\Hom}{\mathrm{Hom}}
\newcommand{\End}{\mathrm{End}}
\newcommand{\Aut}{\mathrm{Aut}}

\newcommand{\Mod}{\mathrm{Mod}}
\newcommand{\Obj}{\mathrm{Obj}\,}
\newcommand{\Rep}{\mathrm{Rep}}
\newcommand{\GLoc}{G\!-\!\mathrm{Loc}}
\newcommand{\Loc}{\mathrm{Loc}}
\newcommand{\Pic}{\mathrm{Pic}}

\def\endexem{\hfill{$\Box$}\medskip}
\newcommand{\qed}{\ \hfill $\blacksquare$\medskip}

\newcommand{\npb}{Nucl. Phys. \1B}
\newcommand{\cmp}{Commun. Math. Phys. }

\newcommand{\rmp}{Rev. Math. Phys. }

\theoremstyle{change}
\theoremheaderfont{\scshape}
\newtheorem{defin}{Definition}[section]
\newtheorem{defprop}[defin]{Definition/Proposition}
\newtheorem{lemma}[defin]{Lemma}
\newtheorem{prop}[defin]{Proposition}
\newtheorem{theorem}[defin]{Theorem}
\newtheorem{coro}[defin]{Corollary}
\newtheorem{conj}[defin]{Conjecture}
\theorembodyfont{\rmfamily}
\newtheorem{rema}[defin]{Remark}
\newtheorem{noname}[defin]{}

\newcommand{\bdefin}{\begin{defin}}
\newcommand{\bdefprop}{\begin{defprop}}
\newcommand{\blemma}{\begin{lemma}}
\newcommand{\bprop}{\begin{prop}}
\newcommand{\btheor}{\begin{theorem}}
\newcommand{\bcoro}{\begin{coro}}
\newcommand{\edefin}{\end{defin}}
\newcommand{\edefprop}{\end{defprop}}
\newcommand{\elemma}{\end{lemma}}
\newcommand{\eprop}{\end{prop}}
\newcommand{\etheor}{\end{theorem}}
\newcommand{\ecoro}{\end{coro}}
\newcommand{\bconj}{\begin{conj}}
\newcommand{\econj}{\end{conj}}
\newcommand{\brem}{\begin{rema}}
\newcommand{\erem}{\endexem\end{rema}}
\newcommand{\bnix}{\begin{noname}}
\newcommand{\enix}{\end{noname}}

\newcommand{\prf}{\noindent{\it Proof. }}

\def\mobj#1{\raise .4\unitlens\hbox{\put(0,0){$#1$}}}

\def\mychi{\raise 2pt\hbox{$\chi$}}

\begin{document}
\maketitle\noindent

\numberwithin{equation}{section}

\center{{\it Dedicated to Detlev Buchholz on the occasion of his sixtieth birthday}} \\ 

\abstract{The aim of the paper is twofold. First, we show that a quantum field
theory $A$ living on the line and having a group $G$ of inner symmetries gives rise to a
category $\GLoc A$ of twisted representations. This category is a braided crossed
$G$-category in the sense of Turaev \cite{t2}. Its degree zero subcategory is braided and
equivalent to the usual representation category $\Rep\,A$. Combining this with \cite{klm},
where $\Rep\,A$ was proven to be modular for a nice class of rational conformal models,
and with the construction of invariants of G-manifolds in \cite{t2}, we obtain an
equivariant version of the following chain of constructions: Rational CFT $\leadsto$
modular category $\leadsto$ 3-manifold invariant.

Secondly, we study the relation between $\GLoc A$ and the braided (in the usual sense) 
representation category $\Rep\,A^G$ of the orbifold theory $A^G$. We prove the equivalence
$\Rep\,A^G\simeq(\GLoc A)^G$, which is a rigorous implementation of the insight that one
needs to take the twisted representations of $A$ into account in order to determine
$\Rep\,A^G$. In the opposite direction we have $\GLoc A\simeq\Rep\,A^G\rtimes\2S$, where
$\2S\subset\Rep\,A^G$ is the full subcategory of representations of $A^G$ contained in
the vacuum representation of $A$, and $\rtimes$ refers to the Galois extensions of braided
tensor categories of \cite{mue06,mue13}. 

Under the assumptions that $A$ is completely rational and $G$ is finite we prove that $A$
has $g$-twisted representations for every $g\in G$ and that the sum over the squared 
dimensions of the simple $g$-twisted representations for fixed $g$ equals $\dim\Rep\,A$.
In the holomorphic case (where $\Rep\,A\simeq\mathrm{Vect}_\7C$) this allows to
classify the possible categories $\GLoc A$ and to clarify the r\^ole of the twisted
quantum doubles $D^\omega(G)$ in this context, as will be done in a sequel. We conclude
with some remarks on non-holomorphic orbifolds and surprising counterexamples concerning
permutation orbifolds.}


\section{Introduction}
It is generally accepted that a chiral conformal field theory (CFT) should have a braided
tensor category of representations, cf.\ e.g.\ \cite{ms}. In order to turn this idea into 
rigorous mathematics one needs an axiomatic formulation of chiral CFTs and their
representations, the most popular framework presently being the one of vertex operator
algebras (VOAs), cf.\ \cite{kac}. It is, however, quite difficult to define a tensor
product of representations of a VOA, let alone to construct a braiding. These difficulties
do not arise in the operator algebraic approach to CFT, reviewed e.g.\ in
\cite{gafr}. (For the general setting see \cite{haag}.) In
the latter approach it has even been possible to give a model-independent proof of
modularity (in the sense of \cite{t1}) of the representation category for a natural class
of rational CFTs \cite{klm}. This class contains the $SU(n)$ WZW-models and the Virasoro
models for $c<1$ and it is closed w.r.t.\ direct products, finite extensions and
subtheories and coset constructions. Knowing modularity of $\Rep\,A$ for rational chiral
CFTs is very satisfactory, since it provides a rigorous way of associating an invariant of
3-manifolds with the latter \cite{t1}. 

It should be mentioned that the strengths and weaknesses of the two axiomatic approaches
are somewhat complementary. The operator algebraic approach has failed so far to reproduce
all the insights concerning the conformal characters afforded by other approaches. (A
promising step towards a fusion of the two axiomatic approaches has been taken in 
\cite{wass}.)

Given a quantum field theory (QFT) $A$, conformal or not, it is interesting to consider
actions of a group $G$ by global symmetries, i.e.\ by automorphisms commuting with the
space-time symmetry. In this situation it is natural to study the relation between the
categories $\Rep\,A$ and $\Rep\,A^G$, where $A^G$ is the $G$-fixed subtheory of $A$. In
view of the connection with string theory, in which the fixpoint theory has a geometric
interpretation, one usually speaks of `orbifold theories'.

In fact, for a quantum field theory $A$ in Minkowski space of $d\ge 2+1$ dimensions
and a certain category $DHR(A)$ of representations \cite{dhr3} -- admittedly too small to
be physically realistic -- the following have been shown \cite{dr2}: (1) $DHR(A)$ is
symmetric monoidal, semisimple and rigid,  (2) there exists a compact group $G$ such that 
$DHR(A)\simeq\Rep\,G$, (3) there exists a QFT $F$ on which $G$ acts by global symmetries
and such that (4) $F^G\cong A$, (5) the vacuum representation of $F$, restricted to $A$,
contains all irreducible representations in $DHR(A)$, (6) all intermediate theories 
$A\subset B\subset F$ are of the form $B=F^H$ for some closed $H\lhd G$, and (7) $DHR(F)$
is trivial. All this should be understood as a Galois theory for quantum fields. 

These results cannot possibly hold in low-dimensional CFT for the simple reason that a
non-trivial modular category is never symmetric. Turning to models with symmetry group
$G$, we will see that $G$ acts on the category $\Rep\,A$ and that $\Rep\,A^G$ contains the
$G$-fixed subcategory $(\Rep\,A)^G$ as a full subcategory. (The objects of the latter are
precisely the representations of $A^G$ that are contained in the restriction to $A^G$ of a
representation of $A$.) Now it is known from models, cf.\ e.g.\ \cite{dvvv}, that
$(\Rep\,A)^G\not\simeq\Rep\,A^G$ whenever $G$ is non-trivial. This can be quantified as 
$\dim\Rep\,A^G=|G|\dim(\Rep\,A)^G=|G|^2\dim\Rep\,A$, cf.\ e.g.\
\cite{x2,mue08}. Furthermore, it 
has been known at least since \cite{dvvv} that $\Rep\,A^G$ is not determined completely by
$\Rep\,A$. This is true  even in the simplest case, where $\Rep\,A$ is trivial but
$\Rep\,A^G$ depends on an additional piece of information pertaining to the `twisted
representations' of $A$. (Traditionally, cf.\ in particular \cite{dvvv,dw,dpr}, it is
believed that this piece of information is an element of $H^3(G,\7T)$, but the situation
is considerably more complicated as we indicate in Subsection \ref{ss-holom} and will be 
elaborated further in a sequel \cite{mue17} to this work.  

Already this simplest case shows that a systematic approach is needed. It turns out that
the right structure to use are the braided crossed G-categories recently introduced for
the purposes of algebraic \cite{carr} and differential \cite{t2} topology. Roughly
speaking, a crossed $G$-category is a tensor category carrying a $G$-grading $\del$ (on
the objects) and a compatible $G$-action $\gamma$. A braiding is a family of isomorphisms 
$(c_{X,Y}: X\otimes Y\rarr {}^XY\otimes X)$, where ${}^XY=\gamma_{\del{X}}(Y)$, satisfying
a suitably generalized form of the braid identities. In Section \ref{s-crossed} we will
show that a QFT on the line carrying a G-action defines a braided crossed G-category
$\GLoc A$ whose degree zero part is $\Rep\,A$. After some further preparation it will turn
out that the additional information contained in $\GLoc A$ is precisely what is needed in
order to compute $\Rep\,A^G$. On the one hand, it is easy to define a `restriction functor'
$R: (\GLoc A)^G\rarr\Rep A^G$, cf.\ Subsection \ref{ss-prelim}. On the other hand, the
procedure of `$\alpha$-induction' from \cite{lre,x0,be123} provides a functor 
$E: \Rep A^G\rarr(\GLoc A)^G$ that is inverse to $R$, proving the braided equivalence
\begin{equation} \Rep\,A^G \simeq (\GLoc A)^G. \label{e-orbif}\end{equation}
Yet more can be said. We recall that given a semisimple rigid braided tensor category
$\2C$ over an algebraically closed field of characteristic zero and a full symmetric
subcategory $\2S$ that is even (all objects have twist $+1$ and thus there exists a
compact group $G$ such that $\2S\simeq\Rep\,G$) there exists a tensor category
$\2C\rtimes\2S$ together with a faithful tensor functor $\iota: \2C\rarr\2C\rtimes\2S$. 
$\2C\rtimes\2S$ is braided if $\2S$ is contained in the center $\2Z_2(\2C)$ of $\2C$
\cite{brug,mue06} and a braided crossed G-category in general \cite{mue13,kir1}. Applying
this to the full subcategory $\2S\subset\Rep\,A^G$ of those representations that are
contained in the vacuum representation of $A$, we show that the functor $E$ factors as  
$E=(\Rep A^G\stackrel{\iota}{\longrightarrow}\Rep
A^G\rtimes\2S\stackrel{F}{\longrightarrow}\GLoc A)$, where $F:\Rep
A^G\rtimes\2S\rarr\GLoc A$ is a full and faithful functor of braided crossed
G-categories. For finite $G$ we prove the latter to be an equivalence: 
\begin{equation} \GLoc A \ \simeq\ \Rep\,A^G\rtimes\2S. \label{e-crossed}\end{equation}
Thus the pair $(\Rep\,A^G,\2S)$ contains the same information as $\GLoc A$
(with its structure as braided crossed G-category).
We conclude that the categorical framework of \cite{mue13} and the quantum field
theoretical setting of Section \ref{s-crossed} are closely related. 

In \cite{klm} it was proven that $\Rep\,A$ is a modular category \cite{t1} if $A$ is
completely rational. In Section \ref{s-cplrtl} we use this result to prove that a
completely rational theory carrying a finite symmetry $G$ always admits $g$-twisted
representations for every $g\in G$. This is an analogue of a similar result \cite{dm} for
vertex operator algebras. (However, two issued must be noted. First, it is not yet known
when a finite orbifold $\2V^G$ of a -- suitably defined -- rational VOA $\2V$ is again
rational, making it at present necessary to assume rationality of $\2V^G$. Secondly, no
full construction of a braided G-crossed category of twisted representations has been
given in the VOA framework.) In fact we have the stronger result  
\[ \sum_{X_i\in(\GLoc A)_g} d(X_i)^2 = \sum_{X_i\in\Loc A} d(X_i)^2
   \ =: \,\dim\Loc A \quad\quad\forall g\in G, \]
where the summations run over the isoclasses of simple objects in the respective categories.

Let us briefly mention some interesting related works. In the operator algebraic
setting, conformal orbifold models were considered in particular in \cite{x2,lx,klx}. 
In \cite{x2} it is shown that $A^G$ is completely rational if $A$ is completely rational
and $G$ is finite, a result that we will use. The other works consider orbifolds in affine
models of CFT, giving a fairly complete analysis of $\Rep\,A^G$. The overlap with our
model independent categorical analysis is small.
Concerning the VOA setting we limit ourselves to mentioning \cite{dm,dy} where suitably
defined twisted representations of $A$ are considered and their existence is proven for all
$g\in G$. Also holomorphic orbifolds are considered. The works \cite{ko,kir1,kir3} are
predominantly concerned with categorical considerations, but the connection with VOAs and
their orbifolds is outlined in \cite[Section 5]{ko}, a more detailed treatment being
announced. \cite[II]{kir1} and \cite{kir3} concern similar matters as \cite{mue06,mue13} 
from a somewhat different perspective. All in all it seems fair to say, however, that no
complete proofs of analogues of our Theorems \ref{theor-M1},\ \ref{theor-M2} and
\ref{theor-M3} for VOAs have been published.

The paper is organized as follows. In Section \ref{s-crossed} we show that a chiral
conformal field theory $A$ carrying a $G$-action gives rise to a braided crossed
$G$-category $\GLoc A$ of (twisted) representations. Even though the construction 
is a straightforward generalization of the procedure in the ungraded case, we give
complete details in order to make the constructions accessible to readers who are
unfamiliar with algebraic QFT. We first consider theories on the line, requiring only the
minimal set of axioms necessary to define $\GLoc A$. We then turn to theories on the
circle, establish the connection between the two settings and review the results of
\cite{klm} on completely rational theories. In Section \ref{s-galois} we study the
relation between the category $\GLoc A$ and the representation category $\Rep\,A^G$ of
the orbifold theory $A^G$, proving (\ref{e-orbif}) and (\ref{e-crossed}). 
In Section \ref{s-cplrtl} we focus on completely rational CFTs \cite{klm} and finite
groups, obtaining stronger results. We give a preliminary discussion of the `holomorphic'
case where $\Rep\,A$ is trivial. A complete analysis of this case is in preparation and
will appear elsewhere \cite{mue17}. We conclude with some comments and counterexamples
concerning  orbifolds of non-holomorphic models.

Most results of this paper were announced in \cite{mue08}, which seems to be the first
reference to point out the relevance of braided crossed G-categories in the context of
orbifold CFT.


\section{Braided Crossed G-Categories in Chiral CFT} \label{s-crossed}
\subsection{QFT on $\7R$ and twisted representations} 
In this subsection we consider QFTs living on the line $\7R$. We begin with some
definitions. Let $\2K$ be the set of intervals in $\7R$, i.e.\ the bounded connected open
subsets of $\7R$. For $I,J\in\2K$ we write $I<J$ and $I>J$ if $I\subset(-\infty,\inf J)$
or $I\subset(\sup J,+\infty)$, respectively. We write $I^\perp=\7R-\ol{I}$. 

For any Hilbert space $\2H$, $\2B(\2H)$ is the set of bounded linear operators on $\2H$,
and for $M\subset\2B(\2H)$ we write $M^*=\{ x^*\ | \ x\in M\}$ and 
$M'=\{ x\in\2B(\2H)\ | \ xy=yx\ \forall y\in M\}$. A {\it von Neumann algebra} (on $\2H$)
is a set $M\subset\2B(\2H)$ such that $M=M^*=M''$, thus in particular it is a unital
$*$-algebra. A {\it factor} is a von Neumann algebra $M$ such that 
$Z(M)\equiv M\cap M'=\7C\11$. A factor $M$ (on a separable Hilbert space) is of {\it type
III} iff for every $p=p^2=p^*\in M$ there exists $v\in M$ such that $v^*v=\11,\ vv^*=p$.
If $M,N$ are von Neumann algebras then $M\vee N$ is the smallest von Neumann algebra
containing $M\cup N$, in fact: $M\vee N=(M'\cap N')'$.

\bdefin \label{def-R}
A QFT on $\7R$ is a triple $(\2H_0,A,\Omega)$, usually simply denoted by $A$, where
\begin{enumerate}
\item $\2H_0$ is a separable Hilbert space with a distinguished non-zero vector $\Omega$,
\item $A$ is an assignment $\2K\ni I\mapsto A(I)\subset\2B(\2H_0)$, where $A(I)$ is a type
III factor.
\end{enumerate}
These data are required to satisfy
\begin{itemize}
\item Isotony: $\ I\subset J\ \ \impl \  A(I)\subset A(J)$,
\item Locality: $I\subset J^\perp\ \impl \  A(I)\subset A(J)'$,
\item Irreducibility: $\vee_{I\in\2K} A(I)=\2B(\2H_0)\ $ (equivalently,
$\cap_{I\in\2K}  A(I)'=\7C\11$),
\item Strong additivity: $A(I)\vee A(J)=A(K)$ whenever $I,J\in\2K$ are adjacent, i.e.\
$\ol{I}\cap\ol{J}=\{p\}$, and $K=I\cup J\cup\{p\}$,
\item Haag duality $A(I^\perp)'=A(I)$ for all $I\in\2K$,
\end{itemize}
where we have used the unital $*$-algebras
\bean A_\infty &=& \bigcup_{I\in\2K} A(I) \ \subset\2B(\2H_0), \\
   A(I^\perp) &=& \mathrm{Alg}\, \{A(J),\ J\in\2K, \ I\cap J=\emptyset \} \ 
  \subset A_\infty. \eean
\edefin

\brem 1. Note that $A_\infty$ is the algebraic inductive limit, no closure is involved.
We have $Z(A_\infty)=\7C\11$ as a consequence of the fact that the $A(I)$ are factors.

2. The above axioms are designed to permit a rapid derivation of the desired categorical
structure. In Subsection \ref{ss-chiral} we will consider a set of axioms that is more
natural from the mathematical as well as physical perspective.
\erem

Our aim is now to associate a strict braided crossed G-category $\GLoc\,A$ to any QFT on
$\7R$ equipped with a G-action on $A$ in the sense of the following

\bdefin \label{def-action}
Let $(\2H_0,A,\Omega)$ be a QFT on $\7R$. A topological group $G$ acts on $A$ if there is
a strongly continuous unitary representation $V: G\rarr\2U(\2H_0)$ such that 
\begin{enumerate}
\item $\beta_g(A(I))=A(I)\ \forall g\in G,\ I\in\2K$, where $\beta_g(x)=V(g)xV(g)^*$.
\item $V(g)\Omega=\Omega$.
\item If $\beta_g\restr A(I)=\id$ for some $I\in\2K$ then $g=e$.
\end{enumerate}
\edefin

\brem 1. Condition 3 will be crucial for the definition of the G-grading on $\GLoc\,A$.

2. In this section the topology of $G$ is not taken into account. In Section
\ref{s-galois} we will mostly be interested in finite groups, but we will also comment on
infinite compact groups. 
\erem

The subsequent considerations are straightforward
generalizations of the well known theory \cite{dhr3, frs1,frs2} for $G=\{e\}$. Since
modifications of the latter are needed throughout -- and also in the interest of the
non-expert reader -- we prefer to develop the case for non-trivial G from scratch.
Readers who are unfamiliar with the following well-known result are encouraged to do the
easy verifications. (We stick to the tradition of denoting the objects of $\End\,B$ by
lower case Greek letters.)

\bdefprop \label{defprop-end}
Let $B$ be a unital $*$-subalgebra of $\2B(\2H)$. Let $\End\,B$ be the category whose
objects $\rho, \sigma,\ldots$ are unital $*$-algebra homomorphisms from $B$ into itself. With
\bean \Hom(\rho,\sigma) &=& \{ s\in B\ | \ s\rho(x)=\sigma(x)s\ \forall x\in B\}, \\
 t\circ s &=& ts, \quad s\in\Hom(\rho,\sigma), t\in\Hom(\sigma,\eta), \\
 \rho\otimes\sigma &=& \rho(\sigma(\cdot)), \\
 s\otimes t &=& s\rho(t)=\rho'(t)s,\quad s\in\Hom(\rho,\rho'),\ t\in\Hom(\sigma,\sigma'), 
\eean
$\End\,B$ is a $\7C$-linear strict tensor category with unit $\11=\id_B$ and positive
$*$-operation. We have $\End\11=Z(B)$.
\edefprop

We now turn to the definition of $\GLoc\,A$ as a full subcategory of $\End\,A_\infty$.

\bdefin 
Let $I\in\2K,\ g\in G$. An object $\rho\in\End\,A_\infty$ is called
$g$-localized in $I$ if 
\[\begin{array}{cccc} \rho(x) &=& x & \quad \forall J<I,\ x\in A(J), \\
   \rho(x) &=& \beta_g(x) & \quad \forall J>I,\ x\in A(J).
\end{array}\]
$\rho$ is $g$-localized if it is $g$-localized in some $I\in\2K$.
A $g$-localized $\rho\in\End\,A_\infty$ is transportable if for every $J\in\2K$
there exists $\rho'\in\End\,A_\infty$, $g$-localized in $J$, such that $\rho\cong\rho'$
(in the sense of unitary equivalence in $\End\,A_\infty$). 
\edefin

\brem \label{rem_act}
1. If $\rho$ is $g$-localized in $I$ and $J\supset I$ then $\rho$ is $g$-localized
in $J$. 

2. Direct sums of transportable morphisms are transportable.

3. If $\rho$ is $g$-localized and $h$-localized then $g=h$. Proof: By 1., there exists 
$I\in\2K$ such that $\rho$ is $g$-localized in $I$ and $h$-localized in $I$. If
$J>I$ then $\rho\restr A(J)=\beta_g=\beta_h$, and condition 3 of Definition
\ref{def-action} implies $g=h$.
\erem

\bdefin \label{def_GRep}
$\GLoc\,A$ is the full subcategory of $\End\,A_\infty$ whose objects are finite
direct sums of $G$-localized transportable objects of $\End\,A_\infty$. Thus
$\rho\in\End\,A_\infty$ is in $\GLoc\,A$ iff there exists a finite set $\Delta$ and,
for all $i\in\Delta$, there exist
$g_i\in G$, $\rho_i\in\End\,A_\infty$ $g_i$-localized transportable,
and $v_i\in\Hom(\rho_i,\rho)$ such that $v_i^*\circ v_j=\delta_{ij}$ and
\[ \rho=\sum_i v_i\,\rho_i(\cdot)\,v_i^*. \]
We say $\rho\in\GLoc\,A$ is $G$-localized in $I\in\2K$ if there exists a
decomposition as above where all $\rho_i$ are $g_i$-localized in $I$ and transportable and  
$v_i\in A(I)\ \forall i$. 

For $g\in G$, let $(\GLoc\,A)_g$ be the full subcategories of $\GLoc\,A$ consisting of 
those $\rho$ that are $g$-localized, and let $(\GLoc\,A)_{\mathrm{hom}}$ be the union of
the $(\GLoc\,A)_g, g\in G$. We write $\Loc\,A=(\GLoc\,A)_e$.

For $g\in G$ define $\gamma_g\in\Aut(\GLoc\,A)$ by 
\bean \gamma_g(\rho) &=& \beta_g\rho\beta_g^{-1}, \\
  \gamma_g(s) &=& \beta_g(s),  \quad s\in\Hom(\rho,\sigma)\subset A_\infty.
\eean
\edefin

\bdefin \label{def_crossed} 
Let $G$ be a (discrete) group. A strict crossed G-category is a strict tensor category
$\2D$ together with 
\begin{itemize}
\item a full tensor subcategory $\2D_{\mathrm{hom}}\subset\2D$ of homogeneous objects,
\item a map $\del: \Obj\,\2D_{\mathrm{hom}}\rarr G$ constant on isomorphism classes,
\item a homomorphism $\gamma: G\rarr\Aut\,\2D$ (monoidal self-isomorphisms of $\2D$)
\end{itemize}
such that 
\begin{enumerate}
\item $\del(X\otimes Y)=\del X\,\del Y$ for all $X,Y\in\2D_{\mathrm{hom}}$. 
\item $\gamma_g(\2D_h)\subset\2D_{ghg^{-1}}$, where $\2D_g\subset\2D_{\mathrm{hom}}$ is
the full subcategory $\del^{-1}(g)$.
\end{enumerate}
If $\2D$ is additive we require that every object of $\2D$ be a direct sum of objects in 
$\2D_{\mathrm{hom}}$. 
\edefin

\bprop $\GLoc\,A$ is a $\7C$-linear crossed G-category with $\End\11=\7C\id_\11$,
positive $*$-operation, direct sums and subobjects (i.e.\ orthogonal projections split).
\eprop

\prf The categories $(\GLoc\,A)_g,\ g\in G$ are mutually disjoint by Remark 
\ref{rem_act}.3. This allows to define the map $\del: \Obj(\GLoc\,A)_{\mathrm{hom}}\rarr G$
required by Definition \ref{def_crossed}. If $\rho$ is $g$-localized in $I$ and $\sigma$
is $h$-localized in $J$ then $\rho\otimes\sigma=\rho\sigma$ is $gh$-localized in any 
$K\in\2K,\ K\supset I\cup J$. Thus $\GLoc\,A$ is a tensor subcategory of $\End\,A_\infty$
and condition 1 of Definition \ref{def_crossed} holds. By construction, $\GLoc\,A$ is
additive and every object is a finite direct sum of homogeneous objects. 
It is obvious that $\gamma_g$ commutes with $\circ$ and with $\otimes$ on objects. Now,
\[ \gamma_g(s\otimes t)=\beta_g(s)\beta_g(\rho(t))
   =\beta_g(s)\gamma_g(\rho)(\beta_g(t))  =\gamma_g(s)\otimes\gamma_g(t). \]
Furthermore, if $s\in\Hom(\rho,\sigma)$ then
$\beta_g(s)\,\beta_g\rho(x)=\beta_g\sigma(x)\,\beta_g(s)$,
and replacing $x\rarr\beta_g^{-1}(x)$ we find
$\beta_g(s)\in\Hom(\gamma_g(\rho),\gamma_g(\sigma))$. Thus $\gamma_g$ it is a strict
monoidal automorphism of $\GLoc\,A$. Obviously, the map $g\mapsto\gamma_g$ is a
homomorphism. If $\rho$ is $h$-localized in $I$ and $J>I$ then 
\[ \gamma_g(\rho)\restr A(J)=\beta_g\rho\beta_g^{-1}=\beta_g\beta_h\beta_g^{-1}, \]
thus $\gamma_g(\rho)$ is $ghg^{-1}$-localized in $I$, thus condition 2 of Definition 
\ref{def_crossed} is verified. 

$\11=\id_{A_\infty}$ is $e$-localized, thus in $\GLoc\,A$ and
$\End\11=Z(A_\infty)=\7C\id_\11$. Let $p=p^2=p^*\in\End(\rho)$. There exists
$I\in\2K$ such that $p\in A(I)$, and by the type III property, cf.\ \ref{conseq}.1,
we find $v\in A(I)$ such that $vv^*=p, v^*v=\11$. Defining $\rho_1=v^*\rho(\cdot)v$ we
have $v\in\Hom(\rho_1,\rho)$, thus $\GLoc\,A$ has subobjects. Finally, for any finite set 
$\Delta$ and any $I\in\2K$ we can find $v_i\in A(I), i\in\Delta$ such that 
$\sum_i v_iv_i^*=\11,\, v_i^*v_j=\delta_{ij}\11$. If $\rho_i\in\GLoc\,A$ we find that
$\rho=\sum_i v_i\rho_i(\cdot)v_i^*$ is a direct sum.  
\qed

\brem Due to the fact that we consider only unital $\rho\in\End A_\infty$, the category
$\GLoc\,A$ does not have zero objects, thus cannot be additive or abelian. This could be
remedied by dropping the unitality condition, but we refrain from doing so since it would
unnecessarily complicate the analysis without any real gains. 
\erem


\subsection{The braiding}
Before we can construct a braiding for $\GLoc\,A$ some preparations are needed. 

\blemma \label{l-br1}
If $\rho$ is $g$-localized in $I$ then $\rho(A(I))\subset A(I)$ and $\rho\restr A(I)$ is
normal. 
\elemma

\prf Let $J<I$ or $J>I$. We have either $\rho\restr A(J)=\id$ or 
$\rho\restr A(J)=\beta_g$. In both cases $\rho(A(J))=A(J)$, implying
$\rho(A(I^\perp))=A(I^\perp)$. Applying $\rho$ to the equation $[A(I),A(I^\perp)]=\{0\}$ 
expressing locality we obtain $[\rho(A(I)),A(I^\perp)]=\{0\}$, or 
$\rho(A(I))\subset A(I^\perp)'=A(I)$, where we appealed to Haag duality on $\7R$. The last 
claim follows from the fact that every unital $*$-endomorphism of a type III factor with
separable predual is automatically normal.
\qed

\blemma \label{l-br0}
Let $\rho,\sigma$ be $g$-localized in $I$. Then $\Hom(\rho,\sigma)\subset A(I)$.
\elemma

\prf Let $s\in\Hom(\rho,\sigma)$. Let $J<I$ and $x\in A(J)$. Then 
$sx=s\rho(x)=\sigma(x)s=xs$, thus $s\in A(J)'$. If $J>I$ and $x\in A(J)$ we find
$s\beta_g(x)=s\rho(x)=\sigma(x)s=\beta_g(x)s$. Since $\beta_g(A(J))=A(J)$ we again have
$s\in A(J)'$. Thus $s\in A(I^\perp)'=A(I)$, by Haag duality on $\7R$.
\qed

\blemma \label{l-br2}
Let $\rho_i\in\GLoc\,A,\ i=1,2$ be $g_i$-localized in $I_i$, where
$I_1<I_2$. Then 
\begin{equation} \label{eq-br3}
 \rho_1\otimes \rho_2=\gamma_{g_1}(\rho_2)\otimes \rho_1. 
\end{equation}
\elemma

\prf We have $I_1=(a,b), I_2=(c,d)$ where $b\le c$. Let $u<a, v>d$ and define 
$K=(u,c), L=(b,v)$. For $x\in A(K)$ we have $\rho_2(x)=x$ and therefore
$\rho_1\rho_2(x)=\rho_1(x)$. By Lemma \ref{l-br1} we have $\rho_1(x)\in A(K)$, and
since $\gamma_{g_1}(\rho_2)$ is $g_1g_2g_1^{-1}$-localized in $I_2$ we find
$\gamma_{g_1}(\rho_2)(\rho_1(x))=\rho_1(x)$. Thus (\ref{eq-br3}) holds for $x\in A(K)$. 
Consider now $x\in A(L)$. By Lemma \ref{l-br1} we have $\rho_2(x)\in A(L)$ and thus
$\rho_1\rho_2(x)=\beta_{g_1}\rho_2(x)$. On the other hand, $\rho_1(x)=\beta_{g_1}(x)$
and therefore
\[ \gamma_{g_1}(\rho_2)\rho_1(x)=\beta_{g_1}\rho_2\beta_{g_1}^{-1}\beta_{g_1}(x)
   =\beta_{g_1}\rho_2(x), \]
thus (\ref{eq-br3}) also holds for $x\in A(L)$. By strong additivity, 
$A(K)\vee A(L)=A(u,v)$, and by local normality of $\rho_1$ and $\rho_2$, (\ref{eq-br3})
holds on $A(u,v)$ whenever $u<a, v>d$, and therefore on all of $A_\infty$. 
\qed

\brem If one drops the assumption of strong additivity then instead of Lemma \ref{l-br1}
one still has $\rho(A(J))\subset A(J)$ for every $J\supset\ol{I}$. Lemma \ref{l-br2} still
holds provided $I_1<I_2$ and $\ol{I_1}\cap\ol{I_2}=\emptyset$.
\erem

Recall that for homogeneous $\sigma$ we write
${}^\sigma\!\!\rho=\gamma_{\del(\sigma)}(\rho)$ as in \cite{t2}.

\bdefin \label{def_braid}
A braiding for a crossed G-category $\2D$ is a family of isomorphisms 
$c_{X,Y}: X\otimes Y\rarr{}^X\!Y\otimes X$, defined for all 
$X\in\2D_{\mathrm{hom}},\ Y\in\2D$, such that
\begin{itemize}
\item[(i)] the diagram
\begin{equation} \begin{diagram}
  X\otimes Y & \rTo^{s\otimes t} & X'\otimes Y' \\
  \dTo_{c_{X,Y}} && \dTo^{c_{X',Y'}} \\
  {}^X\!Y\otimes X & \rTo_{{}^X\!t\otimes s} & {}^{X'}\!Y'\otimes X'
\end{diagram} \label{eq-br0}\end{equation}
commutes for all $s: X\rarr X'$ and $t: Y\rarr Y'$,
\item[(ii)] for all $X,Y\in\2D_{\mathrm{hom}}, \ Z,T\in\2D$ we have
\bea c_{X,Z\otimes T} &=& \id_{{}^X\!Z}\otimes c_{X,T}\mcirc c_{X,Z}\otimes\id_T, 
        \label{eq-br1}\\
   c_{X\otimes Y,Z} &=& c_{X,{}^Y\!Z}\otimes\id_Y\mcirc\id_X\otimes c_{Y,Z}, \label{eq-br2}
\eea
\item[(iii)] for all $X\in\2D_{\mathrm{hom}},\ Y\in\2D$ and $k\in G$ we have
\begin{equation} \gamma_k(c_{X,Y})=c_{\gamma_k(X),\gamma_k(Y)}. \label{eq-br4}\end{equation}
\end{itemize}
\edefin

\bprop \label{lemma-br}
$\GLoc\,A$ admits a unitary braiding $c$. If $\rho_1, \rho_2$ are localized as in Lemma 
\ref{l-br2} then
$c_{\rho_1,\rho_2}=\id_{\rho_1\otimes\rho_2}=\id_{{}^{\rho_1}\!\rho_2\otimes\rho_1}$. 
\eprop

\prf Let $\rho\in(\GLoc\,A)_g$, $\sigma\in\GLoc\,A$ be $G$-localized in $I,J\in\2K$,
respectively. Let $\widetilde{I}<J$. By transportability we can find
$\widetilde{\rho}\in(\GLoc\,A)_g$ localized in $\widetilde{I}$ and a unitary
$u\in\Hom(\rho,\widetilde{\rho})$. By Lemma \ref{l-br2} we have
$\widetilde{\rho}\otimes\sigma=\gamma_{g}(\sigma)\otimes\rho$, thus the composite
\[ \begin{diagram} c_{\rho,\sigma}: & 
\rho\otimes\sigma & \rTo^{u\otimes\id_\sigma}\widetilde{\rho}\otimes\sigma \ &\equiv& 
   \gamma_g(\sigma)\otimes\widetilde{\rho} & \rTo^{\id_{\gamma_g(\sigma)}\otimes u^*} & 
   \gamma_g(\sigma)\otimes\rho 
\end{diagram} \]
is unitary and a candidate for the braiding. As an element of $A_\infty$,
$c_{\rho,\sigma}=\gamma_g(\sigma)(u^*)u=\beta_g\sigma\beta_g^{-1}(u^*)u$.
In order to show that $c_{\rho,\sigma}$ is independent of the choices involved 
pick $\tilde{\tilde{\rho}}\in(\GLoc\,A)_g$ $g$-localized in $\widetilde{I}$ (we may assume
the same localization interval since $\rho$ localized in $\widetilde{I}$ is also localized
in $\tilde{\tilde{I}}\supset\widetilde{I}$) and a unitary
$\widetilde{u}\in\Hom(\rho,\tilde{\tilde{\rho}})$. In view of Lemma \ref{l-br0} we have
$u\widetilde{u}^*\in\Hom(\tilde{\tilde{\rho}},\widetilde{\rho})\subset A(\widetilde{I})$,
implying $\gamma_g(\sigma)(u\widetilde{u}^*)=u\widetilde{u}^*$. The computation
\bean  c_{\rho,\sigma} &=& \gamma_g(\sigma)(u^*)u
   =\gamma_g(\sigma)(u^*)(u\widetilde{u}^*)(\widetilde{u}u^*)u  \\
   &=& \gamma_g(\sigma)(u^*(u\widetilde{u}^*))(\widetilde{u}u^*)u
   =\gamma_g(\sigma)(\widetilde{u}^*)\widetilde{u}
   = \widetilde{c}_{\rho,\sigma}
\eean
shows that $c_{\rho,\sigma}$ is independent of the chosen $\widetilde{\rho}$ and
$u\in\Hom(\rho,\widetilde{\rho})$.

Now consider $\sigma,\sigma'\in\GLoc\,A$ $G$-localized in $J$, $\rho\in(\GLoc\,A)_g$ and 
$t\in\Hom(\sigma,\sigma')$. We pick $\widetilde{I}<J$, $\widetilde{\rho}$ $g$-localized in
$\widetilde{I}$ and a unitary $u\in\Hom(\rho,\widetilde{\rho})$. We define
$c_{\rho,\sigma}=\gamma_g(\sigma)(u^*)u$ and $c_{\rho,\sigma'}=\gamma_g(\sigma')(u^*)u$ as
above. The computation
\[\begin{array}{ccccc}  
      c_{\rho,\sigma'}\mcirc\id_\rho\otimes t
   &=& \gamma_g(\sigma')(u^*)u\,\rho(t) 
   &=& \beta_g\sigma'\beta_g^{-1}(u^*)u\,\rho(t)  \\
   &=& \beta_g\sigma'\beta_g^{-1}(u^*)\widetilde{\rho}(t)u 
   &=& \beta_g\sigma'\beta_g^{-1}(u^*)\beta_g(t)u \\
   &=& \beta_g[\sigma'\beta_g^{-1}(u^*)t]u 
   &=& \beta_g[t\sigma\beta_g^{-1}(u^*)]u \\
   &=& \beta_g(t)\beta_g\sigma\beta_g^{-1}(u^*)u 
   &=& \beta_g(t)\gamma_g(\sigma)(u^*)u \\
   &=& \beta_g(t)\otimes\id_\rho\mcirc c_{\rho,\sigma}
\end{array}\]
proves naturality (\ref{eq-br0}) of $c_{\rho,\sigma}$ w.r.t.\ $\sigma$. (In the fourth
step $\widetilde{\rho}(t)=\beta_g(t)$ is due to $t\in\Hom(\sigma,\sigma')\subset A(J)$, cf.\
Lemma \ref{l-br0}, and the fact that $\rho'$ is $g$-localized in $\widetilde{I}<J$.)

Next, let $\rho,\rho'\in(\GLoc\,A)_g$, $s\in\Hom(\rho,\rho')$ and let $\sigma\in\GLoc\,A$ be
$G$-localized in $J$. Pick $\widetilde{I}<J$, $\widetilde{\rho},\widetilde{\rho}'$
$g$-localized in $\widetilde{I}$ and unitaries
$u\in\Hom(\rho,\widetilde{\rho}),u'\in\Hom(\rho',\widetilde{\rho}')$. Then
\[\begin{array}{ccccc}  
  c_{\rho',\sigma} \mcirc s\otimes\id_\sigma 
  &=& \gamma_g(\sigma)(u'^*)u'\,s
  &=& \gamma_g(\sigma)(u'^*)(u'su^*)u \\
  &=& \gamma_g(\sigma)(u'^*(u'su^*))u
  &=& \gamma_g(\sigma)(su^*)u \\
  &=& \gamma_g(\sigma)(s)\gamma_g(\sigma)(u^*)u
  &=& \id_{\gamma_g(\sigma)}\otimes s\mcirc c_{\rho,\sigma}
\end{array}\]
proves naturality of $c_{\rho,\sigma}$ w.r.t.\ $\rho$. (Here we used the fact that
$\widetilde{\rho},\widetilde{\rho}'$ are $g$-localized in $\widetilde{I}$, implying 
$u'su^*\in\Hom(\widetilde{\rho},\widetilde{\rho}')\subset A(\widetilde{I})$ by Lemma
\ref{l-br0} and finally $\gamma_g(\sigma)(u'su^*)=u'su^*$.)

Next, let $\rho\in(\GLoc\,A)_g$ and let $\sigma,\eta\in\GLoc\,A$ be $G$-localized in $J$.
We pick $\widetilde{\rho}$ $g$-localized in $\widetilde{I}<J$ and a unitary
$u\in\Hom(\rho,\widetilde{\rho})$. Then
\bean c_{\rho,\sigma\otimes\eta} &=& \gamma_g(\sigma\eta)(u^*)u \\
   &=& \gamma_g(\sigma\eta)(u^*)\,\gamma_g(\sigma)(u)\,\gamma_g(\sigma)(u^*)\,u \\
   &=& \gamma_g(\sigma)[\gamma_g(\eta)(u^*)u]\gamma_g(\sigma)(u^*)u \\
   &=& \id_{\gamma_g(\sigma)}\otimes c_{\rho,\eta}\mcirc c_{\rho,\sigma}\otimes\id_\eta \eean
proves the braid relation (\ref{eq-br1}). 

Finally, let $\rho\in(\GLoc\,A)_g, \sigma\in(\GLoc\,A)_h$ and let $\eta\in\GLoc\,A$ be
$G$-localized in $J$. Pick $\widetilde{\rho}\in(\GLoc\,A)_g, \widetilde{\sigma}\in(\GLoc\,A)_h$
$G$-localized in $\widetilde{I}<J$ and unitaries 
$u\in\Hom(\rho,\widetilde{\rho}), v\in\Hom(\sigma,\widetilde{\sigma})$. Then
$w=u\rho(v)=\widetilde{\rho(v)}u\in\Hom(\rho\sigma,\widetilde{\rho}\widetilde{\sigma})$, thus
\bean c_{\rho\otimes\sigma,\eta}
   &=& \gamma_{gh}(\eta)(w^*)w \\
   &=&\gamma_{gh}(\eta)(u^*\widetilde{\rho}(v^*))u\rho(v) \\
   &=& \gamma_{gh}(\eta)(u^*)\gamma_{gh}(\eta)\widetilde{\rho}(v^*)u\rho(v) \\
   &=& \gamma_{gh}(\eta)(u^*)\widetilde{\rho}\gamma_h(\eta)(v^*)u\rho(v) \\
   &=& \gamma_{gh}(\eta)(u^*)u\,\rho[ \gamma_h(\eta)(v^*)v] \\
   &=& \gamma_g(\gamma_h(\eta))(u^*)u\ \rho[ \gamma_h(\eta)(v^*)v]  \\
   &=& c_{\rho,\gamma_h(\eta)}\otimes\id_\sigma\mcirc\id_\rho\otimes c_{\sigma,\eta},
\eean
where we used $\widetilde{\rho}\gamma_h(\eta)=\gamma_{gh}(\eta)\widetilde{\rho}$, cf.\ Lemma
\ref{l-br2}, proves (\ref{eq-br2}). The last claim follows from 
$\rho_1\rho_2={}^{\rho_2}\!\rho_1\rho_2$, cf.\ Lemma \ref{l-br2} and the fact that we may
take $\widetilde{\rho}=\rho$ and $u=\id_\rho$ in the definition of $c_{\rho,\sigma}$.

It remains to show the covariance (\ref{eq-br4}) of the braiding. Recall that
$c_{\rho,\sigma}\in\Hom(\rho\otimes\sigma,\gamma_g(\sigma)\otimes\rho)$ was defined as 
as $\id_{\gamma_g(\sigma)}\otimes u^*\,\circ\,u\otimes\id_\sigma$ for suitable $u$.
Applying the functor $\gamma_k$ we obtain 
\[ \id_{\gamma_{kgk^{-1}}(\gamma_k(\sigma))}\otimes \gamma_k(u)^*\,\circ\,
   \gamma_k(u)\otimes\id_{\gamma_k(\sigma)} \in 
   \Hom(\gamma_k(\rho)\otimes\gamma_k(\sigma),
   \gamma_{kgk^{-1}}(\gamma_k(\sigma))\otimes\gamma_k(\rho)), \]
where $\gamma_k(u)\in\Hom(\gamma_k(\rho),\gamma_k(\widetilde{\rho}))$. Since this is of the
same form as $c_{\gamma_k(\rho),\gamma_k(\sigma)}$ and since the braiding is independent
of the choice of the intertwiner $u$, (\ref{eq-br4}) follows.
\qed


\subsection{Semisimplicity and rigidity}
In view of Lemma \ref{l-br1} we can define

\bdefin $\GLoc_fA$ is the full tensor subcategory of $\GLoc\,A$ of those objects
$\rho$ satisfying $[A(I):\rho(A(I))]<\infty$ whenever $\rho$ is $g$-localized in $I$. 
\edefin

The following is proven by an adaptation of the approach of \cite{gl}.

\bprop $\GLoc_fA$ is semisimple (in the sense that every object is a finite direct sum of
(absolutely) simple objects). Every object of $\GLoc_fA$ has a conjugate in the sense of
\cite{lro} and $\GLoc_fA$ is spherical \cite{bw}.
\eprop

\prf By standard subfactor theory, $[M:\rho(M)]<\infty$ implies that the von Neumann
algebra $M\cap\rho(M)'=\End\,\rho$ is finite dimensional, thus a multi matrix
algebra. This implies semisimplicity since $\GLoc_fA$ has direct sums and subobjects.

Clearly, it is sufficient to show that simple objects have conjugates, thus we consider
$\rho\in(\GLoc_fA)_g$ $g$-localized in $I$.
By the Reeh-Schlieder property \ref{conseq}.3, cf.\ e.g.\ \cite{gafr}, the vacuum $\Omega$
is cyclic and separating for every $A(I), \ I\in\2K$, giving rise to antilinear involutions
$\2J_I=\2J_{(A(I),\Omega)}$ on $\2H_0$, the modular conjugations.  Conditions 1-2 in
Definition \ref{def-action} imply $V(g)\2J_I=\2J_IV(g)$ for all $I\in\2K, g\in G$.
For $z\in\7R$ and $K=(z,\infty)$ it is known \cite{gl,gafr} that 
$j_K: x\mapsto\2J_Kx\2J_K$ maps $A(I)$ onto $A(r_zI)$, where $r_z: \7R\rarr\7R$ is the
reflection about $z$. Thus $j_K$ is an antilinear involutive automorphism of
$A_\infty$. Choosing $z$ to be in the right hand complement of $I$, the geometry is as
follows: 
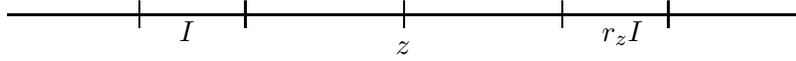
\[ \begin{picture}(300,30)(-150,-15)
\put(0,-5){\line(0,1){10}}
\put(-60,-5){\line(0,1){10}}
\put(-100,-5){\line(0,1){10}}
\put(60,-5){\line(0,1){10}}
\put(100,-5){\line(0,1){10}}
\put(-85,-10){$I$}
\put(75,-10){$r_zI$}
\put(-3,-15){$z$}
\thicklines
\put(-150,0){\line(1,0){300}}
\end{picture}
\]
Let $\widetilde{\rho}$ be $g$-localized in $r_zI$ and $u\in\Hom(\widetilde{\rho},\rho)$ unitary. 
Dropping the subscript $z$ and defining 
\[ \ol{\rho}=j\widetilde{\rho}j\beta_g^{-1}\ \in\End A_\infty \]
it is clear that $\ol{\rho}$ is $g^{-1}$-localized in $I$. It is easy to see that
$d(\ol{\rho})=d(\rho)$ and that $\ol{\rho}$ is transportable, thus in
$(\GLoc_fA)_{g^-1}$. 

Now consider the subalgebras
\[ A_1=\bigcup_{I\in\2K \atop I\subset(-\infty,z)} A(I), \quad\quad
   A_2=\bigcup_{I\in\2K \atop I\subset(z,\infty)} A(I)\]
of $A_\infty$. We have $A_1'=A_2''=\2JA_1''\2J$ . In view of $\widetilde{\rho}\restr A_1=\id$
and $\rho\restr A_2=\beta_g=\mathrm{Ad}\,V(g)$ we have
\bean \rho\restr A_1 &=& u\,\widetilde{\rho}(\cdot)\, u^* = u\,\cdot\, u^*, \\
   \widetilde{\rho}\restr A_2 &=& u^*\,\rho(\cdot)\,u=u^*\,\beta_g(\cdot)\,u
   =u^*V(g)\,\cdot\,V(g)^*u. \eean
We therefore find
\[ \rho\ol{\rho}\restr A_1=\rho j\widetilde{\rho}j\beta_g^{-1}\restr A_1=
   \mathrm{Ad}\,u\2J\,u^*V(g)\,\2JV(g)^* = \mathrm{Ad}\,u\2Ju^*\2J, \]
where we used the commutativity of $\2J$ and $V(g)$. Since the above expressions for 
$\rho\restr A_1, \widetilde{\rho}\restr A_2,\rho\ol{\rho}\restr A_1$ are ultraweakly continuous
they uniquely extend to the weak closures $A_1'', A_2'', A_1''$, respectively. Now,
\bean  u\2Ju^*\, \rho(A_1)''\,u\2Ju^*  &=& u\2J\,\widetilde{\rho}(A_1)''\,\2Ju^* 
  = u\2J\,A_1''\,\2Ju^*  = u\,A_2''\,u^* \\ 
   &=& (u\,A_2'\,u^*)'=(u\,A_1''\,u^*)'=\rho(A_1'')'=\rho(A_1)'.
\eean
Thus, $\widetilde{\2J}=u\2Ju^*$ is an antiunitary involution whose adjoint action maps
$\rho(A_1)''$ onto $\rho(A_1)'$. Furthermore, $u\Omega$ is cyclic and separating for
$\rho(A_1)''$ and we have $(u\2Ju^*)(u\Omega)=u\2J\Omega=u\Omega$ and 
$(\rho(x)\widetilde{\2J}\rho(x)\widetilde{\2J}u\Omega,u\Omega)=(x\2Jx\2J\Omega,\Omega)\ge 0\ \forall x\in A_1''$.
Thus $\widetilde{\2J}$ is \cite[Exercise 9.6.52]{kr2} the modular conjugation corresponding to
the pair $(\rho(A_1)'',u\Omega)$, and therefore 
\[ x\mapsto\rho\ol{\rho}(x)=\2J_{(\rho(A_1)'',u\Omega)}\2J_{(A_1'',\Omega)}\, x\,
   \2J_{(A_1'',\Omega)}\2J_{(\rho(A_1)'',u\Omega)} \]
is a canonical endomorphism $\gamma: A_1''\hookrightarrow\rho(A_1)''$ \cite{lo2}. Since 
$[A_1'':\rho(A_1'')]=[A(I):\rho(A(I))]=d(\rho)^2$ is finite by assumption, $\gamma$
contains \cite{lo2} the 
identity morphism, to wit there is $V\in A_1''$ such that $Vx=\rho\ol{\rho}(x)V$ for all
$x\in A_1''$. Since $\rho\ol{\rho}$ is (e-)localized in $I$, Lemma \ref{l-br0} implies 
$V\in A(I)$, thus the equation $Vx=\rho\ol{\rho}(x)V$ also holds for $x\in A(I')$, and  
strong additivity together with local normality of $\rho, \ol{\rho}$ imply that it holds
for all $x\in A_\infty$. Thus $\11=\id_{A_\infty}\prec\rho\ol{\rho}$, and $\ol{\rho}$ is a
conjugate, in the sense of \cite{lro}, of $\rho$ in the tensor $*$-category
$\GLoc_fA$. Choosing a conjugate or dual $\ol{\rho}$ for every $\rho\in\GLoc_fA$ and
duality morphisms $e:\ol{\rho}\otimes\rho\rarr\11,\ \11\rarr\rho\otimes\ol{\rho}$
satisfying the triangular equations we may consider $\GLoc_fA$ as a spherical category.  
\qed

\brem Every object $\rho$ in a spherical or $C^*$-category with simple unit has a dimension
$d(\rho)$ living in the ground field, $\7C$ in the present situation. This dimension of an
object localized in $I$ is related to the index by the following result of Longo \cite{lo2}:
\[ d(\rho)=[A(I):\rho(A(I))]^{1/2}. \]
\erem

Summarizing the preceding discussion we have 

\btheor  \label{theor-M1}
$\GLoc\,A$ is a braided crossed G-category and $\GLoc_fA$ is a rigid semisimple braided
crossed G-category.
\etheor

\brem
1. It is obvious that for any braided G-crossed category $\2D$, the degree zero subcategory 
$\2D_e$ is a braided tensor category. In the case at hand, $\Loc\,A=(\GLoc\,A)_e$
is the familiar category of transportable localized morphisms defined in \cite{frs1}. 
But for non-trivial symmetries $G$, the category $\GLoc\,A$ contains information that
cannot be obtained from $\Loc\,A$.

2. The closest precedent to our above considerations can be found in \cite{khr2}. There,
however, several restrictive assumptions were made, in particular only abelian groups $G$
were considered. Under these assumptions the $G$-crossed structure essentially trivializes.
\erem


\subsection{Chiral conformal QFT on $S^1$} \label{ss-chiral}
In this subsection we briefly recall the main facts pertinent to chiral conformal field
theories on $S^1$ and their representations, focusing in particular the completely
rational models introduced and analyzed in \cite{klm}. While nothing in this subsection is
new, we include the material since it will be essential in what follows.

Let $\2I$ be the set of intervals in $S^1$, i.e.\ connected open non-empty and non-dense 
subsets of $S^1$. ($\2I$ can be identified with the set 
$\{ (x,y)\in S^1\times S^1\ | \ x\ne y\}$.) For every $J\subset S^1$, $J'$ is the interior
of the complement of $J$. This clearly defines an involution on $\2I$.

\bdefin \label{def-chiral}
A chiral conformal field theory is a quadruple $(\2H_0,A,U,\Omega)$, usually
simply denoted by $A$, where 
\begin{enumerate}
\item $\2H_0$ is a separable Hilbert space with a distinguished non-zero vector $\Omega$,
\item $A$ is an assignment $\2I\ni I\mapsto  A(I)$, where $A(I)$ is a von Neumann algebra
on $\2H_0$.
\item $U$ is a strongly continuous unitary representation of the M\"obius group 
$PSU(1,1)=SU(1,1)/\{\11,-\11\}$, i.e.\ the group of those fractional linear maps
$\7C\rarr\7C$ which map the circle into itself, on $\2H_0$.
\end{enumerate}
These data must satisfy
\begin{itemize}
\item Isotony: $\ I\subset J\ \ \impl \  A(I)\subset A(J)$,
\item Locality: $I\subset J'\ \impl \  A(I)\subset A(J)'$,
\item Irreducibility: $\vee_{I\in\2I} A(I)=\2B(\2H_0)\ $ (equivalently,
$\cap_{I\in\2I}  A(I)'=\7C\11$),
\item Covariance: $U(a) A(I)U(a)^*= A(aI)\quad\forall a\in PSU(1,1),\ I\in\2I$,
\item Positive energy: $L_0\ge 0$, where $L_0$ is the generator of the rotation subgroup
of $PSU(1,1)$,
\item Vacuum: every vector in $\2H_0$ which is invariant under the action of $PSU(1,1)$ is a
multiple of $\Omega$.
\end{itemize}
\edefin

\bnix \label{conseq}
For consequences of these axioms see, e.g., \cite{gafr}. We limit ourselves to listing
some facts: 

\begin{enumerate}
\item Type: The von Neumann algebra $A(I)$ is a factor of type III (in fact III$_1$)
for every $I\in\2I$.
\item Haag duality: $A(I)'= A(I')\ \  \forall I\in\2I$.
\item Reeh-Schlieder property: $\ol{ A(I)\Omega}=\ol{ A(I)'\Omega}=\2H_0 \ \ \forall I\in\2I$.
\item The modular groups and conjugations associated with $( A(I),\Omega)$ have a
geometric meaning, cf.\ \cite{bgl,gafr} for details.
\item Additivity: If $I,J\in\2I$ are such that $I\cap J, I\cup J\in\2I$ then
$A(I)\vee A(J)=A(I\cup J)$. 
\end{enumerate}
\enix

In order to obtain stronger results we introduce two further axioms.

\bdefin Two intervals $I, J\in\2I$ are called adjacent if their closures intersect in
exactly one point. A chiral CFT satisfies strong additivity if
\[ I, J\  \mbox{adjacent} \quad\impl\quad  A(I)\vee A(J)= A(\ol{I\cup J}^0). \]
A chiral CFT satisfies the split property if the map
\[ m:\  A(I)\otimes_{alg} A(J)\rarr A(I)\vee A(J),\quad x\otimes y\mapsto xy \]
extends to an isomorphism of von Neumann algebras whenever $I, J\in\2I$ satisfy
$\ol{I}\cap\ol{J}=\emptyset$.
\edefin

\brem By M\"obius covariance strong additivity holds in general if it holds for one pair
$I, J$ of adjacent intervals. Furthermore, every CFT can be extended canonically to one
satisfying strong additivity. The split property is implied by the property
$Tr e^{-\beta L_0}<\infty \ \forall\beta>0$. The latter property and strong additivity
have been verified in all known rational models.
\erem

\bdefin \label{def_rep}
A representation $\pi$ of $A$ on a Hilbert space $\2H$ is a family $\{\pi_I, I\in\2I\}$,
where $\pi_I$ is a unital $*$-representation of $A(I)$ on $\2H$ such that  
\begin{equation} \label{eq-rep}
   I\subset J \quad \impl \quad \pi_J\restr A(I)=\pi_I. 
\end{equation}
$\pi$ is called covariant if there is a positive energy representation $U_\pi$ of the 
universal covering group $\widehat{PSU(1,1)}$ of the M\"obius group on $\2H$ such that  
\[ U_\pi(a)\pi_I(x)U_\pi(a)^* = \pi_{aI}(U(a)xU(a)^*) \quad \forall
   a\in\widehat{PSU(1,1)}, \ I\in\2I. \]
We denote by $\Rep\,A$ the $C^*$-category of all representations on separable Hilbert
spaces, with bounded intertwiners as morphisms.
\edefin

\bdefprop \label{def_rep2}
If $A$ satisfies strong additivity and $\pi$ is a representation then the Jones index 
of the inclusion $\pi_I(A(I))\subset\pi_{I'}(A(I'))$ does not depend on $I\in\2I$ and we
define the dimension
\[ d(\pi)=[\pi_{I'}(A(I')):\pi_I(A(I))]^{1/2} \ \in[1,\infty]. \]
We define $\Rep_f A$ to be the the full subcategory of $\Rep\,A$ of those representations 
satisfying $d(\pi)<\infty$.
\edefprop

As just defined, $\Rep\,A$ and $\Rep_fA$ are just $C^*$-categories. In order to obtain the
well known result \cite{frs1,gafr} that the category of all (separable) representations
can be equipped with braided monoidal structure, we need the following: 

\bprop \label{prop-SR}
Every chiral CFT $(\2H_0,A,U,\Omega)$ satisfying strong additivity gives rise to a QFT on 
$\7R$.
\eprop

\prf We arbitrarily pick a point $\infty\in S^1$ and consider
\[ \2I_\infty=\{ I\in\2I\ | \ \infty\not\in\ol{I}\} \]
Identifying $S^1-\{\infty\}$ with $\7R$ by stereographic projection
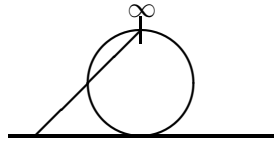
\[
\begin{picture}(300,60)(-150,-30)\thicklines
\put(0,0){\circle{50}}
\put(-50,-20){\line(1,0){100}}
\put(0,15){\line(0,1){10}}
\put(-5,25){$\infty$}
\put(0,20){\line(-1,-1){40}}
\end{picture}
\]
we have a bijection between $\2I_\infty$ and $\2K$. The family $A(I), I\in\2K$ is just the
restriction of $A(I), I\in\2I$ to $I\in\2I_\infty\equiv\2K$. By \ref{conseq}, $A$
satisfies Haag duality on $S^1$, and together with strong additivity (on $S^1$) this
implies Haag duality (on $\7R$) and strong additivity in the sense of Definition
\ref{def-R}.
\qed

\brem The definition of $G$-actions on a chiral CFT on $S^1$ is analogous to Definition
\ref{def-action}, condition 1 now being required for all $I\in\2I$. Conditions 1-2 imply
$V(g)U(a)=U(a)V(g)\ \forall g\in G,\ a\in PSU(1,1)$. (To see this observe that 1-2 imply
that $V(g)$ commutes with the modular groups associated with the pairs $(A(I),\Omega)$ for
any $I\in\2I$. By \ref{conseq}.4 the latter are one-parameter subgroups of $U(PSU(1,1))$
which generate $U(PSU(1,1))$.) Condition 3 now is equivalent to the more convenient axiom
\begin{itemize}
\item[3'.] If $U(g)\in\7C\11$ then $g=e$.
\end{itemize}
(Proof: If $U(g)\in\7C\11$ then $\alpha_g=\id$, thus $g=e$ by 3. Conversely, if $\alpha_g$
acts trivially on some $A(I)$ then $U(g)$ commutes with $A(I)$ and in fact with all $A(I)$
by $V(g)U(a)=U(a)V(g)$. Thus the irreducibility axiom implies $U(g)\in\7C\11$.)
\erem

Given a CFT on $S^1$ and ignoring a possibly present $G$-action we have the categories 
$\Rep\,A$ ($\Rep_fA$) as well as the braided tensor categories $\Loc\,A$ ($\Loc_fA$) 
associated with the restriction of $A$ to $\7R$. The following result, cf.\
\cite[Appendix]{klm}, connects these categories. 

\btheor 
Let $(\2H_0,A,U,\Omega)$ be a chiral CFT satisfying strong additivity. Then there are
equivalences of $*$-categories 
\bean \Loc\,A &\simeq& \Rep\,A, \\
  \Loc_fA &\simeq& \Rep_fA, \eean
where $\Rep_{(f)}A$ refers to the chiral CFT and Definition \ref{def_rep}, whereas
$\Loc_{(f)}A$ refers to the QFT on $\7R$ obtained by restriction and Definition
\ref{def_GRep}. 
\etheor

\prf The strategy is to construct a functor $Q: \Loc\,A\rarr\Rep\,A$ of
$*$-categories and to prove that it is fully faithful and essentially surjective. Let 
$\rho\in\Loc\,A$ be localized in $I\in\2K\equiv\2I_\infty$. Our aim is to define a
representation $\pi=(\pi_I,\ I\in\2I)$ on the Hilbert space $\2H_0$. For every
$J\in\2I_\infty$ we define $\pi_J=\rho\restr A(J)$, considered as a representation on
$\2H_0$. If $\infty\in\ol{J}$ we pick an
interval $K\in\2I_\infty$, $K\cap J=\emptyset$. By transportability of $\rho$ there
exists $\rho'$ localized in $K$ and a unitary $u\in\Hom(\rho,\rho')$. Defining 
$\pi_J=u^*\cdot u$ we need to show that $\pi_J$ is independent of the choices involved. 
Thus let $\rho''$ be localized in $K$ (this may be assumed by making $K$ large enough) and 
$v\in\Hom(\rho,\rho'')$, giving rise to $\pi_J'=v^*\cdot v$. Now, 
$u\circ v^*\in\Hom(\rho'',\rho')$, thus 
$uv^*\in A(K)$ by Lemma \ref{l-br0}, and therefore 
\[ \pi_J'(x)=v^*xv=v^*(vu^*uv^*)xv=v^*vu^*xuv^*v =u^*xu=\pi_J(x), \]
since $x\in A(J)\subset A(K)'$. Having defined $\pi_J$ for all $J\in\2I$ we need to show 
(\ref{eq-rep}) for all $I,J\in\2I$. There are three cases of inclusions $I\subset J$ to be
considered: (i) $I, J\in\2I_\infty$, (ii) $I\in\2I_\infty, J\not\in\2I_\infty$,
(iii) $I,J\not\in\2I_\infty$. Case (i) is trivial since $\pi_I=\pi_J=\rho$, restricted to
$A(I), A(J)$ respectively. Case (iii) is treated by using $K\subset J'$ for the definition
of both $\pi_I, \pi_J$ and appealing to the uniqueness of the latter. In case (ii) we have
$\pi_J=u^*\cdot u$ with $u\in\Hom(\rho,\rho')$, $\rho'$ localized in $K\subset J'$. For
$x\in A(I)$ we have $\pi_J(x)=u^*xu=u^*\rho'(x)u=\rho(x)=\pi_I(x)$, as desired.
This completes the proof of $\pi=\{\pi_J\}\in\Rep\,A$.

Let $\rho_1,\rho_2\in\Loc\,A$ and let $\pi_1,\pi_2$ be the corresponding
representations. We claim that $s\in\Hom(\rho_1,\rho_2)$ implies $s\in\Hom(\pi_1,\pi_2)$. 
Let $\infty\in\ol{J}$, $K\in\2I_\infty, K\cap J=\emptyset$, $\rho'_i$ localized in $K$ and
$u_i\in\Hom(\rho_i,\rho_i')$ unitaries, such that then $\pi_{J,i}=u_i^*\cdot u_i$. We have
$u_2su_1^*\in\Hom(\rho_1',\rho_2')$. Since both $\rho_1', \rho_2'$ are localized in $K$ we
have $u_2su_1^*\in A(K)\subset A(J)'$. Now the computation
\[ s\pi_{J,1}(x)=su_1^*xu_1= u_2^*(u_2su_1^*)xu_1= u_2^*x(u_2su_1^*)u_1
   = u_2^*xu_2s= \pi_{J,2}s \]
shows $s\in\Hom(\pi_{J,1},\pi_{J,2})$. Since this works for all $J$ such that
$\infty\in\ol{J}$ we have $s\in\Hom(\pi_1,\pi_2)$, and we have defined a faithful functor
$Q: \Loc\,A\rarr\Rep\,A$. Obviously, $Q$ is faithful. In view of
$\rho=\pi\restr A_\infty$ it is clear that $s\in\Hom(\pi,\pi')$ implies
$s\in\Hom(\rho,\rho')$, thus $Q$ is full. 

Let now $\pi\in\Rep\,A$ and $I\in\2I$. Then
$\pi_I$ is a unital $*$-representation of $A(I)$ on a separable Hilbert space. Since
$A(I)$ is of type III and $\2H_0$ is separable, $\pi_I$ is unitarily implemented. I.e.\
there exists a unitary $u: \2H_0\rarr\2H_\pi$ such that $\pi_I(x)=uxu^*$ for all 
$x\in A(I)$. Then $(\pi_J')=(u^*\pi_J(\cdot)u)$ is a representation on $\2H_0$ that
satisfies $\pi'\cong\pi$ and $\pi_I'=\pi_{I,0}=\id$. Haag duality (on $S^1$) implies 
$\pi_J(A(J))\subset A(J)$ whenever $J\supset I'$. If we choose $I$ such that $\infty\in I$
then $\pi_J, J\supset I'$ defines an endomorphism $\rho$ of $A_\infty$ whose extension to
a representation $Q(\rho)$ coincides with $\pi'$. Thus $Q$ is essentially surjective and
therefore an equivalence $\Loc\,A\simeq\Rep\,A$. 

Now, $\rho\in\Loc\,A$ is in $\Loc_fA$ iff $d(\rho)=[A(I):\rho(A(I))]^{1/2}<\infty$ 
whenever $\rho$ is localized in $I$. On the other hand, $\pi\in\Rep\,A$ is in $\Rep_fA$ iff
$d(\pi)=[\pi_{I'}(A(I')):\pi_I(A(I))]^{1/2}<\infty$. In view of the above construction it
is clear that $d(\pi)=d(\rho)$ if $\pi$ is the representation corresponding to $\rho$. 
Thus $Q$ restricts to an equivalence $\Loc_fA\simeq\Rep_fA$. 
\qed

Using the equivalence $Q$ the braided monoidal structure of $\Loc_{(f)}A$ can be
transported to $\Rep_{(f)}A$:

\bcoro $\Rep\,A$ ($\Rep_fA$) can be equipped with a (rigid) braided monoidal structure
such that there are equivalences
\bean \Loc\,A &\simeq& \Rep\,A, \\
  \Loc_fA &\simeq& \Rep_fA \eean
of braided monoidal categories. 
\ecoro

\brem \label{rem-3}
1. It is quite obvious that the braided tensor structure on $\Rep\,A$ provided by the 
above constructions is independent, up to equivalence, of the choice of the point
$\infty\in S^1$. For an approach to the representation theory of QFTs on $S^1$ that does
not rely on cutting the circle see \cite{frs2}. The latter, however, seems less suited
for the analysis of $\GLoc\,A$ for non-trivial $G$ since the $g$-localized endomorphisms
of $A_\infty$ do not extend to endomorphisms of the global algebra $\2A_{\mathrm{univ}}$ of
\cite{frs2} if $g\ne e$.

2. Given a chiral CFT $A$, the category $\Rep\,A$ is a very natural object to consider. 
Thus the significance of the degree zero category $(\GLoc\,A)_e$ is plainly evident: It
enables us to endow $\Rep\,A$ with a braided monoidal structure in a considerably easier
way than any known alternative.

3. By contrast, the rest of the category $\GLoc\,A$ has no immediate physical
interpretation. After all, the objects of $(\GLoc_fA)_g$ with $g\ne e$ do not
represent proper representations of $A$ since they `behave discontinuously at $\infty$'.
In fact, it is not difficult to prove that, given two adjacent intervals $I,J\in\2I$ and
$g\ne e$, there exists no representation $\pi$ of $A$ such that $\pi\restr A(I)=\id$ and
$\pi\restr A(J)=\beta_g$. Thus $\rho$, considered as a representation of $A_\infty$,
cannot be extended to a representation of $A$. The main {\it physical} relevance of
$\GLoc_fA$ is that -- in contradistinction to $\Rep_fA$ -- it contains sufficient
information to compute $\Rep_fA^G$. This will be discussed in the next section.

4. On the purely {\it mathematical} side, the category $\GLoc\,A$ may be used to define an
invariant of three dimensional G-manifolds \cite{t2}, i.e.\ 3-manifolds equipped with a
principal G-bundle. As mentioned in the introduction, this provides an equivariant version
of the construction of a 3-manifold invariant from a rational CFT.
\erem

As is well known, there are models, like the $U(1)$ current algebra, that satisfy the
standard axioms including strong additivity and the split property and that have
infinitely many inequivalent irreducible representations. Since in this work we are mainly
interested in rational CFTs we need another axiom to single out the latter.

\bdefprop \cite{klm} Let $A$ satisfy strong additivity and the split property.
Let $I,J\in\2I$ satisfy $\ol{I}\cap\ol{J}=\emptyset$ and write $E=I\cup J$. Then
the index of the inclusion $A(E)\subset A(E')'$ does not depend on $I,J$ and we define
\[ \mu(A)=[A(E')':A(E)]\ \in[1,\infty]. \]
A chiral CFT on $S^1$ is completely rational if it satisfies (a) strong additivity, (b)
the split property and (c) $\mu(A)<\infty$.
\edefprop

\brem 1. Thus every CFT satisfying strong additivity and the split property comes along
with a numerical invariant $\mu(A)\in [1,\infty]$. The models where the latter is finite 
-- the completely rational ones -- are among the best behaved (non-trivial) quantum field
theories, in that very strong results on both their structure and representation theory
have been proven in \cite{klm}. In particular the invariant $\mu(A)$ has a nice
interpretation.

2. All known classes of rational CFTs are completely rational in the above sense. For the 
WZW models connected to loop groups this is proven in \cite{wass,x1}. More importantly,
the class of completely rational models is stable under tensor products and finite
extensions and subtheories, cf.\ Section \ref{s-galois} for more details. This has
applications to orbifold and coset models. 
\erem

\btheor \cite{klm} \label{KLM_main} 
Let $ A$ be a completely rational CFT. Then
\begin{itemize}
\item Every representation of $A$ on a separable Hilbert space is completely reducible,
i.e.\ a direct sum of irreducible representations. (For non-separable representations
this is also true if one assumes local normality, which is automatic in the separable
case, or equivalently covariance.) 
\item Every irreducible separable representation has finite dimension $d(\pi)$, thus
$\Rep_f A$ is just the category of finite direct sums of irreducible representations.
\item The number of unitary equivalence classes of separable irreducible representations
is finite and
\[ \dim\Rep_fA=\mu(A), \]
where $\dim\Rep_fA$ is the sum of the squared dimensions of the simple objects.
\item The braiding of $\Loc_fA\simeq\Rep_fA$ is non-degenerate, thus $\Rep_f A$ is a
unitary modular category in the sense of Turaev \cite{t1}.
\end{itemize}
\etheor


\section{Orbifold Theories and Galois Extensions} \label{s-galois}
\subsection{The restriction functor $R: (\GLoc A)^G\rarr\Loc A^G$} \label{ss-prelim}
After the interlude of the preceding subsection we now return to QFTs defined on $\7R$
with symmetry $G$. (Typically they will be obtained from chiral CFTs on $S^1$ by 
restriction, but in the first subsections this will not be assumed.) Our aim is to
elucidate the relationship between the categories $\GLoc\,A$ and $\Loc\,A^G$, where $A^G$
is the `orbifold' subtheory of $G$-fixpoints in the theory $A$.

\bdefin Let $(\2H,A,\Omega)$ be a QFT on $\7R$ with an action (in the sense of
Definition \ref{def-action}) of a compact group $G$. Let $\2H_0^G$ and $A(I)^G$ be
the fixpoints under the G-action on $\2H_0$ and $A(I)$, respectively. Then the orbifold
theory $A^G$ is the triple $(\2H_0^G,A^G,\Omega)$, where $A^G(I)=A(I)^G\restr\2H_0^G$.
\edefin

\brem \label{rem-orb}
1. The definition relies on $\Omega\in\2H_0^G$ and $A(I)^G\2H_0^G\subset\2H_0^G$ for
all $I\in\2K$. Denoting by $p$ the projector onto $\2H_0^G$, we have
$A^G(I)=A(I)^G\restr\2H_0^G=pA(I)p$, where the right hand side is understood as an algebra
acting on $p\2H_0=\2H_0^G$. Furthermore, since $A(I)^G$ acts faithfully on $\2H_0^G$ we
have algebra isomorphisms $A(I)^G\cong A^G(I)$.

2. It is obvious that the triple $(\2H_0^G,A^G,\Omega)$ satisfies isotony and locality. 
Irreducibility follows by $\vee_{I\in\2K}A^G(I)=p(\vee_{I\in\2K}A(I))p$ together with
$\vee_I A(I)=\2B(\2H_0)$. However, strong additivity and Haag duality of the fixpoint
theory are not automatic. For the time being we will postulate these properties to
hold. Later on we will restrict to settings where this is automatically the case.
\erem

\bnix \label{nix-S}
For later purposes we recall a well known fact about compact group actions on QFTs in the 
present setting. Namely, for every $I\in\2K$, the $G$-action on $A(I)$ has full
$\widehat{G}$-spectrum, \cite{dhr1}. This means that for every isomorphism class 
$\alpha\in\widehat{G}$ of irreducible representations of $G$ there exists a finite
dimensional $G$-stable subspace $V_\alpha\subset A(I)$ on which the $G$-action restricts
to the irrep $\pi_\alpha$. $V_\alpha$ can be taken to be a space of isometries of support
$1$. (This means that $V_\alpha$ admits a basis $\{v_\alpha^i,\ i=1,\ldots,d_\alpha\}$
such that $\sum_i v_\alpha^i{v_\alpha^i}^*=\11$ and
${v_\alpha^i}^*v_\alpha^j=\delta_{ij}\11$.) Furthermore, $A(I)$ is generated by $A(I)^G$
and the spaces $V_\alpha,\alpha\in\widehat{G}$.

These observations have an important consequence for the representation categories of
fixpoint theories \cite{dhr1}. Namely the category $\Loc_fA^G$ contains a full symmetric
subcategory $\2S$ equivalent to the category $\Rep_fG$ of finite dimensional continuous
unitary representations of $G$. The objects in $\2S$ are given by the localized
endomorphisms of $A^G_\infty$ of the form 
$\rho_\alpha(\cdot)=\sum_i v_\alpha^i\cdot{v_\alpha^i}^*$, where $\{v_\alpha^i\}$ is a
space of isometries with support $\11$ in $A(I)$ transforming under the irrep
$\alpha\in\widehat{G}$. (Equivalently, a simple object of $\rho\in\Loc_fA^G$ is in $\2S$
iff the corresponding representation $\pi_0\circ\rho$ of $A^G$ is contained in the
restriction to $A^G$ of the defining (or vacuum) representation of $A$.)
\enix

\bnix \label{rem-orb2}
We now begin our study of the relationship between $\GLoc A$ and $\Loc A^G$. Let 
$(\GLoc A)^G$ denote the $G$-invariant objects and morphisms of $\GLoc A$. By definition
of the $G$-action on $\GLoc A$, $\rho\in(\GLoc A)^G$ implies 
$\rho\circ\beta_g=\beta_g\circ\rho$ for all $g\in G$, thus
$\rho(A^G_\infty)\subset A^G_\infty$. Every $\rho\in\GLoc A$ is $G$-localized in some
interval $I$. In view of Definition \ref{def_GRep} it is obvious that the restriction
$\rho\restr A^G_\infty$ acts trivially on $A(J)$ not only if $J<I$, but also if $J>I$.
Thus $\rho\restr A^G_\infty$ is a localized endomorphism of $A^G_\infty$. Furthermore, if
$\rho,\sigma\in(\GLoc A)^G$ and $s\in\Hom_{(\GLoc A)^G}(\rho,\sigma)$ it is easy to see
that $s\in\Hom_{\Loc A^G}(\rho\restr A^G_\infty,\sigma\restr A^G_\infty)$. This suggests 
that $\rho\restr A^G_\infty\in\Loc A^G$. However, this also requires showing that the
restricted morphism $\rho\restr A^G_\infty$ is transportable {\it by morphisms in $\Loc
A^G$}. This requires some work.
\enix

\bprop \label{prop-restr}
Let $\rho\in(\GLoc A)^G$. Then $\rho\restr A^G_\infty\in\Loc A^G$.
\eprop

\prf By definition, $\rho$ is $G$-localized in some interval $I$. As we have seen in
\ref{rem-orb2}, $\rho\restr A^G_\infty$ is localized in $I$, and it remains to show that 
$\rho\restr A^G_\infty$ is transportable. Let thus $J$ be another interval. By
transportability of $\rho\in\GLoc A$, there exists $\widetilde{\rho}$ that is
$G$-localized in $J$ and a unitary $u\in\Hom_{\GLoc A}(\rho,\widetilde{\rho})$. Define
$\widetilde{\rho}_g=\gamma_g(\widetilde{\rho})=\beta_g\circ\widetilde{\rho}\circ\beta_g^{-1}$.
Since $\gamma_g$ is an automorphism of $\GLoc A$ and $\rho$ is $G$-invariant we have
$\gamma_g(u):=\beta_g(u)\in\Hom_{\GLoc A}(\rho,\widetilde{\rho}_g)$. Defining
$v_g=\beta_g(u)u^*$ we have 
\[ v_{gh}=\beta_{gh}(u)u^*=\beta_g(v_h)\beta_g(u)u^*=\beta_g(v_h)v_g \quad\forall g,h. \] 
Furthermore, $v_g\in\Hom(\widetilde{\rho},\widetilde{\rho}_g)$, and since all
$\widetilde{\rho}_g$ are $G$-localized in $J$, Lemma \ref{l-br0} implies $v_g\in A(J)$. 
Thus $g\mapsto v_g$ is a (strongly continuous) 1-cocycle in $A(J)$. Since $A(J)$ is a type
III factor and the $G$-action has full $\widehat{G}$-spectrum, there exists \cite{suth} a
unitary $w\in A(J)$ such that $v_g=\beta_g(w)w^*$ for all $g\in G$. Defining
$\widehat{\rho}=\mathrm{Ad}\,w^*\circ\widetilde{\rho}$, we have
$w^*u\in\Hom(\rho,\widehat{\rho})$. Now, $\beta_g(u)u^*=\beta_g(w)w^*$ is equivalent to
$\beta_g(w^*u)=w^*u$, thus $w^*u$ is $G$-invariant. Together with the obvious fact that 
$\widehat{\rho}$ is $G$-localized in $J$, this implies $\rho\restr A^G_\infty\in\Loc A^G$. 
\qed

\bcoro Restriction to $A^G_\infty$ provides a strict tensor functor 
$R: (\GLoc A)^G\rarr\Loc A^G$ which is faithful on objects and morphisms.
\ecoro

\prf With the exception of faithfulness, which follows from the isomorphisms 
$A(I)^G\cong A^G(I)$, this is just a restatement of our previous results.
\qed

\brem 1. In Subsection \ref{ss-equiv1} we will show that $R$, when restricted to
$(\GLoc_fA)^G$, is also surjective on morphisms (thus full) and objects. Thus $R$ will  
establish an isomorphism $(\GLoc_f A)^G\cong \Loc_f A^G$. 

2. We comment on our definition of the fixpoint category $\2C^G$ of a category $\2C$ under
a $G$-action. In the literature, cf.\ \cite{tam,kir1,kir3}, one can find a different
notion of fixpoint category, which we denote by $\2C_G$ for the present purposes. Its
objects are pairs $(X,\{u_g, g\in G\})$, where $X$ is an object of $\2C$ and the 
$u_g\in\Hom_\2C(X,\gamma_g(X))$ are isomorphisms making the left diagram in Figure 1
commute. The morphisms between $(X,\{u_g, g\in G\})$ and $(Y,\{v_g, g\in G\})$
are those $s\in\Hom_\2C(X,Y)$ for which the right diagram in Figure 1 commutes. (According
to J.~Bernstein, $\2C_G$ should rather be called the category of $G$-modules in $\2C$.)
\begin{figure}[h] 
\[
\begin{diagram} X & \rTo^{u_g} & \gamma_g(X) \\ & \rdTo_{u_{gh}} & \dTo_{\gamma_g(u_h)} \\ &&
  \gamma_{gh}(X) \end{diagram} 
\quad\quad\quad\quad\quad\quad
\begin{diagram}
X & \rTo^{s} & Y \\ \dTo^{u_g} && \dTo_{v_g} \\ \gamma_g(X) & \rTo_{\gamma_g(s)} &
\gamma_g(Y) 
\end{diagram}
\]
\caption{Objects and Morphisms of $\2C_G$}
\end{figure}
It is clear that $\2C^G$ can be identified with a full subcategory of $\2C_G$ via
$X\mapsto (X,\{ \id\})$, but in general this inclusion need not be an equivalence.
However, it is an equivalence in the case of $\2C=\GLoc A$. To see this, let 
$(\rho,\{u_g\})\in(\GLoc A)_G$. Assume $\rho$ is $G$-localized in $I$. By definition of
$(\GLoc A)_G$, $g\mapsto u_g$ is a 1-cocycle in $A(I)$, and by the above discussion there
exists $w\in A(I)$ such that $u_g=\beta_g(w)w^*$ for all $g\in G$. Defining
$\widetilde{\rho}=\mathrm{Ad}\,w^*\circ\rho$, an easy computation shows
$\widetilde{\rho}\in(\GLoc A)^G$. Since $w:\widetilde{\rho}\rarr\rho$ is an isomorphism,
the inclusion $(\GLoc A)^G\hookrightarrow(\GLoc A)_G$ is essentially surjective, thus an
equivalence. 
\erem


\subsection{The extension functor $E: \Loc A^G\rarr(\GLoc A)^G$} \label{ss-ext}
In view of Remark \ref{rem-orb} we are in a setting where both $A=(\2H_0,A(\cdot),\Omega)$
and $A^G=(\2H_0^G,A^G(\cdot),\Omega)$ are QFTs on $\7R$. In this situation it is well known
that there exists a monoidal functor $E: \Loc\,A^G\rarr\End\,A_\infty$ from the tensor
category of localized transportable endomorphisms of the subtheory $A^G$ to the (not
a priori localized) endomorphisms of the algebra $A_\infty$. There are essentially three
ways to construct such a functor. First, Roberts' method of localized cocycles, cf.\
e.g. \cite{rob}, which is applicable under the weakest set of assumptions. (Neither
finiteness of the extension nor factoriality or Haag duality are required.) Unfortunately,
in this approach it is relatively difficult to make concrete computations, cf.\ however
\cite{cdr}. Secondly, the subfactor approach of Longo and Rehren \cite{lre} as further
studied by Xu, B\"ockenhauer and Evans, cf.\ e.g.\ \cite{x0,be123}. This approach requires
factoriality of the local algebras and finiteness of the extension, but otherwise is very
powerful. Thirdly, there is the approach of \cite{mue04}, which assumes neither
factoriality nor finiteness, but which is restricted to extensions of the form 
$A^G\subset A$. For the present purposes, this is of course no problem.

\btheor \cite{mue04} \label{theor-E}
Let $A=(\2H_0,A(\cdot),\Omega)$ be a QFT on $\7R$ with $G$-action
such that $A^G=(\2H_0^G,A^G(\cdot),\Omega)$ is a QFT on $\7R$. There is a functor 
$E: \Loc\,A^G\rarr\End\,A_\infty$ with the following properties: 
\begin{enumerate}
\item For every $\rho\in\Loc\,A^G$ we have that $E(\rho)$ commutes with the $G$-action
$\beta$, i.e.\ $E(\rho)\in(\End\,A_\infty)^G$. The restriction $E(\rho)\restr A^G_\infty$
coincides with $\rho$. On the arrows, $E$ is the inclusion 
$A^G_\infty\hookrightarrow A_\infty$. Thus $E$ is faithful and injective on the objects.
\item $E$ is strict monoidal. (Recall that $\Loc\,A^G$ and $\End\,A_\infty$ are strict.)
\item If $\rho$ is localized in the interval $I\in\2K$ then $E(\rho)$ is localized in the
half-line $(\inf I,+\infty)$. This requirement makes $E(\rho)$ unique.
\end{enumerate}
\etheor

\noindent{\it Remarks on the proof:} Fix an interval $I\in\2K$. By \ref{nix-S}, we can
find a family $\{V_\alpha\subset A(I), \alpha\in\widehat{G}\}$ of finite dimensional
subspaces of isometries of support $1$ on which the $G$-action restricts to the
irreducible representation $\alpha\in\widehat{G}$. Now the algebra $A(I)$ is generated 
by $A(I)^G$ and the family $\{V_\alpha,\alpha\in\widehat{G}\}$, and $A_\infty$ is
generated by $A^G_\infty$ and the family $\{V_\alpha,\alpha\in\widehat{G}\}$. Furthermore,
$\sigma_\alpha=\sum_i v_\alpha^i\cdot{v_\alpha^i}^*$ is a transportable endomorphism of
$A^G_\infty$ localized in $I$, thus $\sigma_\alpha\in\Loc_fA^G$. Now $E(\rho)$ is
determined by Rehren's prescription \cite{khr3}:
\[ E(\rho)(x)=\left\{ \begin{array}{ll} \rho(x) & x\in A^G_\infty \\ 
   c(\sigma_\alpha,\rho)x & x\in V_\alpha, \end{array} \right. \]
where $c(\sigma_\alpha,\rho)$ is the braiding of the category $\Loc\,A^G$. (The proof of
existence and uniqueness of $E(\rho)$ is given in \cite{mue04}, generalizing the
automorphism case treated in \cite{dr1}. Note that despite the appearances this definition
of $E$ does not depend on the chosen spaces $V_\alpha$.) An the arrows
$\Hom_{\Loc\,A^G}(\rho,\sigma)\subset A^G_\infty$ we define $E$ via the inclusion
$A^G_\infty\hookrightarrow A_\infty$. For the verification of all claimed properties see
\cite[Proposition 3.11]{mue04}.  
\qed

\brem 1. The definition of $E$ does not require $d(\rho)<\infty$. But from now on we will
restrict $E$ to the full subcategory $\Loc_fA^G\subset\Loc_fA^G$.

2. The extension functor $E$ is faithful but not full. Our aim will be to compute
$\Hom_{\End\,A_\infty}(E(\rho),E(\sigma))$, but this will require some categorical
preparations. 
\erem


\subsection{Recollections on Galois extensions of braided tensor categories}\label{ss-galois} 
From the discussion in \ref{nix-S} it is clear that the extension
$E(\rho)\in\End\,A_\infty$ is trivial, i.e.\ isomorphic to a direct sum of
$\dim(\rho)\in\7N$ copies of the tensor unit $\11$, for every $\rho$ in the full symmetric
subcategory $\2S$. It is therefore natural to ask for the universal faithful tensor
functor $\iota: \2C\rarr\2D$ that trivializes a full symmetric subcategory $\2S$ of a
rigid braided tensor category $\2C$. Such a functor has been constructed independently in
\cite{mue06} (without explicit discussion of the universal property) and in \cite{brug}. 
(The motivation of both works was to construct a modular category from a non-modular
braided category by getting rid of the central/degenerate/transparent objects.) A
universal functor $\iota:\2C\rarr\2D$ trivializing $\2S$ exists provided every object in
$\2S$ has trivial twist $\theta(X)$, both approaches relying on the fact \cite{dr6,del}
that under this condition $\2S$ is equivalent to the representation category of a group
$G$, which is finite if $\2S$ is finite and otherwise compact \cite{dr6} or proalgebraic
\cite{del}. In the subsequent discussion we will use the approach of \cite{mue06} since it
was set up with the present application in mind, but we will phrase it in the more
conceptual way expounded in \cite{mue13}.

Given a rigid symmetric tensor $*$-category $\2S$ with simple unit and trivial twists, the
main result of \cite{dr6} tells us that there is a compact group $G$ such that
$\2S\simeq\Rep_fG$. (In our application to the subcategory $\2S\subset\Loc_fA^G$ for an
orbifold CFT $A^G$ we don't need to appeal to the reconstruction theorem since the
equivalence $\2S\simeq\Rep_fG$ is proven already in \cite{dhr1}.) Assuming $\2S$ (and thus
$G$) to be finite we know that there is a commutative strongly separable Frobenius algebra
$(\gamma,m,\eta,\Delta,\ve)$ in $\2S$, where $\gamma$ corresponds to the left regular
representation of $G$ under the equivalence. See \cite{mue09} for the precise definition
and proofs. (More generally, this holds for any finite dimensional semisimple and
cosemisimple Hopf algebra $H$ \cite{mue09}. For infinite compact groups and infinite
dimensional discrete quantum groups one still has an algebra structure $(\gamma,m,\eta)$,
cf.\ \cite{MT}.) The group $G$ can be recovered from the monoid structure
$(\gamma,m,\eta)$ as
\[ G\cong \{ s\in\End \gamma \ | \ s\circ m=m\circ s\otimes s,\ s\circ\eta=\eta\}. \]
Now we define \cite{mue13} a category $\2C\rtimes_0\2S$ with the same
objects and same tensor product of objects as $\2C$, but larger hom-sets:
\[ \Hom_{\2C\rtimes_0\2S}(\rho,\sigma)=\Hom_\2C(\gamma\otimes\rho,\sigma). \]
The compositions $\circ,\otimes$ of morphisms are defined using the Frobenius algebra
structure on $\gamma$. Finally, $\2C\rtimes\2S$ is defined as the idempotent completion
(or Karoubian envelope) of $\2C\rtimes_0\2S$. The latter contains $\2C\rtimes_0\2S$ as a 
full subcategory and is unique up to equivalence, but there also is a well known canonical
model for it. I.e., the objects of $\2C\rtimes\2S$ are pairs $(\rho,p)$, where
$\rho\in\2C\rtimes_0\2S$ and $p=p^2=p^*\in\End_{\2C\rtimes_0\2S}(\rho)$. The morphisms are
given by 
\[ \Hom_{\2C\rtimes\2S}((\rho,p),(\sigma,q))=q\circ\Hom_{\2C\rtimes_0\2S}(\rho,\sigma)\circ q
   =\{ s\in\Hom_{\2C\rtimes_0\2S}(\rho,\sigma) \ | \ s=q\circ s\circ p\}. \]
The inclusion functor $\iota: \2C\rarr\2C\rtimes\2S,\ \rho\mapsto(\rho,\id_\rho)$ has the
desired trivialization property since $\dim\Hom_{\2C\rtimes\2S}(\11,\iota(\rho))=d(\rho)$
for all $\rho\in\2S$. The group $G$ acts on a morphism
$s\in\Hom_{\2C\rtimes\2S}((\rho,p),(\sigma,q))\subset\Hom_\2C(\gamma\otimes\rho,\sigma)$
via $\gamma_g(s)=s\circ g^{-1}\otimes\id_\rho$, where $g\in\Aut(\gamma,m,\eta)\cong G$.
The $G$-fixed subcategory $(\2C\rtimes\2S)^G$ is just the idempotent completion of $\2C$
and thus equivalent to $\2C$. The braiding $c$ of $\2C$ lifts to a braiding of
$\2C\rtimes\2S$ iff all objects of $\2S$ are central, i.e.\
$c(\rho,\sigma)c(\sigma,\rho)=\id$ for all $\rho\in\2S$ and $\sigma\in\2C$. This, however,
will not be the case in the application to QFT. As shown in \cite{mue13}, in the general
case $\2C\rtimes\2S$ is a braided crossed $G$-category. We need one concrete formula from
\cite{mue13}. Namely, if
$p\in\End_{\2C\rtimes\2S}(\rho)\cong\Hom_\2C(\gamma\otimes\rho,\rho)$ is such that 
$(\rho,p)\in\2C\rtimes\2S$ is simple, then the morphism
\begin{equation} \label{e3}
\del(\rho,p)=\quad
\left( \ \ \
\begin{tangle}
\hh\step\coev\\
\hh\frabox{p}\step\id\\
\obj{\eta}\counit\step\hev\\
\end{tangle}
\ \ \right)^{-1}
\cdot\quad
\begin{tangle}
\object{\gamma}\\
\hh\id\Step\coev\\
\x\step\id\obj{\ol{\rho}}\\
\x\step\id\\
\hh\id\step\frabox{p}\step\id\\
\hh\step[-.5]\obj{\Delta}\hstep\hcu\hstep\obj{\rho}\hstep\hev\\
\hh\hstep\id\\
\hstep\object{\gamma}
\end{tangle}
\end{equation}
is an automorphism of the monoid $(\gamma,m,\eta)$, thus an element of $G$. We note for
later use that the numerical factor $(\cdots)^{-1}$ is $d(\rho,p)^{-1}$ and that replacing
the braidings by their duals  
($\begin{tangle}\x\end{tangle}\ \leftrightarrow\ \begin{tangle}\xx\end{tangle}$) gives
the inverse group element.

If the category $\2S$, equivalently the group $G$ are infinite, the above definition of
$\2C\rtimes\2S$ needs to be reconsidered since, e.g., the proof of semisimplicity must be
modified. The original construction of $\2C\rtimes\2S$ in \cite{mue06} does just
that. Using the decomposition $\gamma\cong\oplus_{i\in\widehat{G}}d(\gamma_i)\gamma_i$ of
the regular representation one defines
\begin{equation} \Hom_{\2C\rtimes_0\2S}(\rho,\sigma)=\bigoplus_{i\in\widehat{G}}  
  \Hom_\2C(\gamma_i\otimes\rho,\sigma)\otimes \2H_i, \label{eq-hom}\end{equation}
where $F: \2S\rarr\Rep_fG$ is an equivalence, $\gamma_i\in\2S$ is such that
$F(\gamma_i)\cong\pi_i$ and $H_i$ is the representation space of the irreducible
representation $\pi_i$ of $G$. (It is easily seen that
$\Hom_{\2C\rtimes_0\2S}(\rho,\sigma)$ is finite dimensional for all $\rho,\sigma\in\2C$.) 
Now the compositions $\circ,\otimes$ of morphisms are defined by the formulae
\[ s\boxtimes\psi_k\circ t\boxtimes\psi_l = \bigoplus_{m\in\widehat{G}} \sum_{\alpha=1}^{N_{kl}^m} 
   s\,\circ \,\id_{\gamma_k}\otimes t \,\circ\,
   w_{kl}^{m\alpha}\otimes\id_\rho \, \boxtimes \, K(w_{kl}^{m\alpha})^* 
   (\psi_k\otimes\psi_l), \]
\[  u\boxtimes\psi_k\,\otimes\,w\boxtimes\psi_l = \bigoplus_{m\in\widehat{G}} 
   \sum_{\alpha=1}^{N_{kl}^m} u \otimes v \,\circ\, \id_{\gamma_k} \otimes
   \ve(\gamma_l,\rho_1)\otimes\id_{\rho_2} \,\circ\, 
   w_{kl}^{m\alpha}\otimes\id_{\rho_1\rho_2} \, \boxtimes \,
  K(w_{kl}^{m\alpha})^* (\psi_k\otimes\psi_l), \]
where $k,l\in\widehat{G},\ \psi_k\in\2H_k,\ \psi_l\in\2H_l$,
$t\in\Hom(\gamma_l\otimes\rho,\sigma), s\in\Hom(\gamma_k\otimes\sigma,\delta)$ and 
$u\in\Hom(\gamma_k\otimes\rho_1,\sigma_1), V\in\Hom(\gamma_l\otimes\rho_2,\sigma_2)$. For
for further details and the definition of the $*$-involution, which we don't need here, we
refer to \cite{mue06}. For finite $G$ it is readily verified that the two definitions of
$\2C\rtimes\2S$ given above produce isomorphic categories. If $\Gamma$ is central in
$\2C$, equivalently $c(\rho,\sigma)c(\sigma,\rho)=\id$ for all $\rho\in\2S,\sigma\in\2C$,
then $\2C\rtimes\2S$ inherits the braiding of $\2C$, cf.\ \cite{mue06}. If this is not the
case, $\Gamma-\Mod$ is only a braided crossed G-category \cite{mue13}. 

Before we return to our quantum field theoretic considerations we briefly comment on the
approach of \cite{brug} and the related works \cite{par,ko,kir1,kir3}. As before, one
starts from the (Frobenius) algebra in $\2S$ corresponding to the left regular
representation of $G$. One now considers the category $\Gamma-\Mod$ of left modules over
this algebra. As already observed in \cite{par}, this is a tensor category. Again, if
$\Gamma$ is central in $\2C$ then $\Gamma-\Mod$ is braided \cite{brug}, whereas in general
$\Gamma-\Mod$ is a braided crossed G-category \cite{kir1,kir3}. (The braided degree zero
subcategory coincides with the dyslexic modules of \cite{par}.) In \cite{mue13} an
equivalence of $\2C\rtimes\2S$ and $\Gamma-\Mod$ is proven. In the present investigations
it is more convenient to work with $\2C\rtimes\2S$ since it is strict if $\2C$ is.


\subsection{The isomorphism $\Loc_fA^G \protect\cong (\GLoc_fA)^G$} \label{ss-equiv1}
In Subsection \ref{ss-ext}, the extension functor $E$ was defined on the entire category
$\Loc\,A^G$. It is faithful but not full, and our aim is to obtain a better understanding
of $\Hom_{A_\infty}(E(\rho),E(\sigma))$. 
From now on we will restrict it to the full subcategory $\Loc_fA^G$ of finite
dimensional (thus rigid) objects, and we abbreviate $\2C=\Loc_fA^G$ throughout. 
Furthermore, $\2S\subset\2C$ will denote the full subcategory discussed in \ref{nix-S}. We
recall that $\2S\simeq\Rep_fG$ as symmetric tensor category. Since the definition of
$\2C\rtimes\2S$ in \cite{mue06} was motivated by the formulae \cite{khr3,mue04} for the
intertwiner spaces $\Hom_{A_\infty}(E(\rho),E(\sigma))$, the following is essentially obvious:

\bprop \label{prop-F}
Under the same assumptions on $A$ and $A^G$ and notation as above, the functor
$E:\2C\rarr(\End\,A_\infty)^G$ factors through the canonical inclusion functor 
$\iota: \2C\hookrightarrow\2C\rtimes\2S$, i.e.\ there is a tensor functor 
$F: \2C\rtimes\2S\rarr\End\,A_\infty$ such that 
\[ \begin{diagram} \2C & \rTo^{\iota} & \2C\rtimes \2S \\ & \rdTo_{E} & \dTo_{F} \\ 
   && \End\,A_\infty \end{diagram}\]
commutes. (Note that $F(\2C\rtimes\2S)\not\subset(\End\,A_\infty)^G$.)
The functors 
\[\begin{array}{lccc} E: & \2C & \rarr & (\End\, A_\infty)^G, \\ 
  F: & \2C\rtimes\2S & \rarr & \End\,A_\infty 
\end{array}\]
are faithful and full.
\eprop

\prf First, we define $F$ on the tensor category $\2C\rtimes_0\2S$ of \cite{mue06,mue13},
which has the same objects as $\2C$ but larger hom-sets. We clearly have to put
$F(\rho):=E(\rho)$. Now fix an interval $I\in\2K$ and subspaces $H_i\subset A(I)$ of
isometries on which $G$ acts according to the irrep $\pi_i$. Let $\gamma_i$ be the
endomorphism of $A^G_\infty$ implemented by $H_i$. As stated in \cite{khr3} and proven in
\cite{mue04}, the intertwiner spaces between extensions $E(\rho),E(\sigma)$ is given by
\[ \Hom_{A_\infty}(E(\rho),E(\sigma))=\mathrm{span}_{i\in\widehat{G}}
   \Hom_\2C(\gamma_i\rho,\sigma)H_i\ \subset\ A_\infty. \]
On the one hand, this shows that every $G$-invariant morphism
$s\in\Hom_{\GLoc_fA}(E(\rho),E(\sigma))$ is in $\Hom_{\Loc_fA^G}(\rho,\sigma)$, implying
that $E: \2C\rarr(\End\,A_\infty)^G$ is full. On the other hand, it is clear that these
spaces can be identified with those in the second definition 
(\ref{eq-hom}) of $\Hom_{\2C\rtimes_0\2S}(\rho,\sigma)$. Under this identification, the
compositions $\circ,\otimes$ of morphisms in $\2C\rtimes_0\2S$ go into those in the
category $\End\,A_\infty$ as given in Definition \ref{defprop-end}, as is readily
verified. Thus we have a full and faithful strict tensor functor 
$F_0:\2C\rtimes_0\2S\rarr\End\,A_\infty$ such that $F_0\circ\iota=E$. Now, $\2C\rtimes\2S$
is defined as the completion of $\2C\rtimes_0\2S$ with splitting idempotents. Since the
category $\End\,A_\infty$ has splitting idempotents, the functor $F_0$ extends to a tensor
functor $F:\2C\rtimes\2S\rarr\End\,A_\infty$, uniquely up to natural isomorphism of
functors. However, we give a more concrete prescription. Let $(\rho,p)$ be an object of
$\2C\rtimes\2S$, i.e.\ $\rho\in\Loc_fA$ and $p=p^2=p^*\in\End_{\2C\rtimes_0\2S}(\rho)$. 
Let $I\subset\2K$ be an interval in which $\rho\in\Loc_fA$ is localized. Then Haag
duality implies $p\in A(I)$. Since $A(I)$ is a type III factor (with separable predual) we
can pick $v\in A(I)$ such that $vv^*=p$ and $v^*v=1$. Now we define
$F((\rho,p))(\cdot)=v^*F(\rho)(\cdot)v\in\End\,A_\infty$. This is an algebra 
endomorphism of $A_\infty$ since $vv^*=p\in\Hom_{A_\infty}(E(\rho),E(\rho))$. With this
definition, the functor $F:\2C\rtimes\2S\rarr\End\,A_\infty$ is strongly (but not
strictly) monoidal. 
\qed

In \cite{mue13} it was shown that $\2C\rtimes\2S$ is a braided crossed $G$-category. In
view of the results of Section 2 it is natural to expect that the functor $F$ actually
takes its image in $\GLoc A$ and is a functor of braided crossed $G$-categories. In fact:

\bprop \label{prop-gloc}
Let $A=(\2H_0,A(\cdot),\Omega)$ be as before and $G$ finite. Then 

(i) for every $\rho\in\Loc_fA^G$ we have $E(\rho)\in\GLoc_fA$, thus the extension
$E(\rho)$ is a finite direct sum of endomorphisms $\eta_i$ of $A_\infty$ that act as
symmetries $\beta_{g_i}$ on a half line $[a,+\infty)$. 

(ii) $F(\2C\rtimes\2S)\subset\GLoc_fA$ and $F:\2C\rtimes\2S\rarr\GLoc_fA$ is a
functor of $G$-graded categories, i.e.\ $F((\2C\rtimes\2S)_g)\subset(\GLoc_fA)_g$ for all
$g\in G$. 
\eprop

\prf Claim (i) clearly follows from (ii). In order to prove the latter it is sufficient to
show for every irreducible object $(\rho,p)\in\2C\rtimes\2S$ that
$E((\rho,p))\in\End\,A_\infty$ is $\del(\rho,p)$-localized. Let thus
$\rho\in\2C=\Loc_fA^G$ be localized in the interval $I\in\2K$ and let $p=p^2=p^*$ be a
minimal projection in $\End_{\2C\rtimes_0\2S}(\rho)$. Recall that $F((\rho,p))$ is defined
as $v^*E(\rho)(\cdot)v$, where $v\in A_\infty$ satisfies $vv^*=F(p),\, v^*v=1$. We may
assume that $v\in A(I)$. Let $J\in\2K$ such that $I<J$ and let $H_\gamma\subset A(J)$ be a
subspace of isometries transforming under the left regular representation of $G$. (I.e.,
we have isometries $v_g\in A(I), g\in G$ such 
that $\beta_k(v_g)=v_{kg},\ \sum_g v_vv_g^*=1, v_g^*v_h=\delta_{g,h}1$.) Let 
$\gamma(\cdot)=\sum_g v_g \cdot v_g^*\in\End\,A^G_\infty$ the localized endomorphism
implemented by $H_\gamma$. Thus $H_\gamma=\Hom_A(\11,E(\gamma))$.
Now by Theorem \ref{theor-E} we have, for $x\in H_\gamma$,
\[ F((\rho,p))(x)=v^* c(\gamma,\rho)xv = [ E(\gamma)(v^*)c(\rho,\gamma)c(\gamma,\rho)
  E(\gamma)(v) ] \,x, \] 
where we have used (i) $xv=vx$ (since $x,v$ are localized in the disjoint intervals $I,J$,
respectively), (ii) $c(\rho,\gamma)=1$ (follows by Lemma \ref{lemma-br} since the
localization region of $\rho$ is in the left complement of the localization region of
$\gamma$) and (iii) $E(\gamma)(v)=v$ (since $v\in A(I)$, on which $E(\rho)$ acts
trivially). This expression defines an element of $\Hom_A(F((\rho,p),F(\gamma)F((\rho,p))$.
If $v_1,\ldots,v_{|G|}\in\Hom_A(\11,E(\gamma))$ are such that 
$\sum_i v_iv_i^*=\11,\ v_i^*v_j=\delta_{i,j}\11$ then
\[ v_i^*\,[ E(\gamma)(v^*)c(\rho,\gamma)c(\gamma,\rho) E(\gamma)(v) ] \,x
  \in\End_A(F(\rho,p)). \]
By irreducibility of $F((\rho,p))$ this expression is a multiple of $\id_{F(\rho,p)}$,
thus
\bean d((\rho,p))\,F((\rho,p))(x) &= & 
  d(\rho,p)\,[ E(\gamma)(v^*)c(\rho,\gamma)c(\gamma,\rho) E(\gamma)(v) ] \,x \\
  &=& \sum_i v_i \,
  Tr_{(\rho,p)}(v_i^*\,[ E(\gamma)(v^*)c(\rho,\gamma)c(\gamma,\rho) E(\gamma)(v) ] \,x) \\
&=&
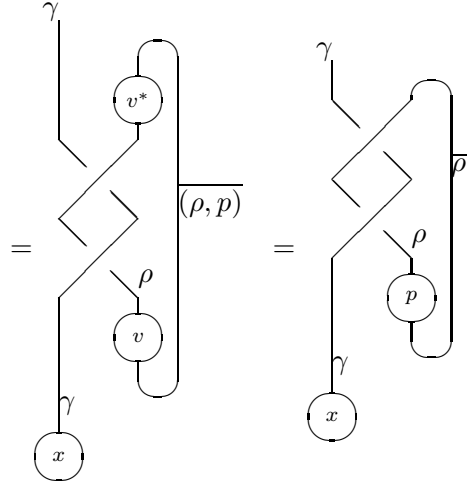
\begin{tangle}
\step[-.2]\object{\gamma}\\
\hh\id\Step\hcoev\\
\id\Step\O{v^*}\step\id\\
\xx\step\id\obj{\overline{(\rho,p)}}\\
\xx\mobj{\rho}\step\id\\
\id\Step\O{v}\step\id\\
\hh\obj{\gamma}\id\Step\hev\\
\QQ{x}
\end{tangle}
\quad\quad\quad=\quad\begin{tangle}
\step[-.2]\object{\gamma}\\
\hh\id\Step\hcoev\\
\xx\step\id\obj{\overline{\rho}}\\\\
\xx\mobj{\rho}\step\id\\
\id\Step\O{p}\step\id\\
\hh\obj{\gamma}\id\Step\hev\\
\QQ{x}
\end{tangle}
\eean
Now we express this as a diagram in $\2C$ in terms of the representers
$x\in\Hom_\2C(\gamma,\gamma)$ and $p\in\Hom_\2C(\gamma\otimes\rho,\rho)$. By definition of
$\2C\rtimes\2S$ we obtain
\[
d((\rho,p))\,F((\rho,p))(x)=\quad 
\begin{tangle}
\step[-.2]\object{\gamma}\\
\hh\id\Step\hcoev\\
\xx\step\id\\
\xx\mobj{\rho}\step\id\\
\hh\id\step\frabox{p}\step\id\\
\hxx\step\hev\\
\id\step\O{x}\hh\\
\hh\step[-.5]\obj{\Delta}\hstep\hcu\\
\hstep\object{\gamma}
\end{tangle}
\quad=\quad
\begin{tangle}
\step[-.2]\object{\gamma}\\
\O{x}\Step\hcoev\\
\xx\step\id\\
\xx\mobj{\rho}\step\id\\
\hh\id\step\frabox{p}\step\id\\
\hh\hcu\step\hev\\
\hstep\object{\gamma}
\end{tangle}
\]
where we have used the commutativity $\Delta=c(\gamma,\gamma)\circ\Delta$.
Thus by (\ref{e3}) and \cite{mue13} we have  
\[ F((\rho,p))(x)= x\circ \del((\rho,p))^{-1}, \]
where $\del((\rho,p))\in\Aut(\gamma,m,\eta)$ is the degree of $(\rho,p)$. Recalling that
the action of $g\in\Aut(\gamma,m,\eta)$ on the morphism
$s\in\Hom_\2C(\gamma\otimes\rho,\sigma)\cong\Hom_{\2C\rtimes\2S}(\rho,\sigma)$ was defined
as $\gamma_g(s)=s\circ g^{-1}\otimes\id_\rho$, we see that 
$F((\rho,p))(x)=\gamma_{\del(\rho,p)}(x)$. Thus $F((\rho,p))\in\End\,A_\infty$ is
$\del(\rho,p)$-localized in the sense of Section \ref{s-crossed}, as
claimed. Transportability of $E((\rho,p))$ follows from transportability of $\rho$. Thus
$E((\rho,p))\in\GLoc_fA$, and the same clearly follows for the non-simple objects of
$\Loc_fA$. The above computations have also shown that the functor $F$ respects the
G-gradings of $\2C\rtimes\2S$ and $\GLoc\,A$ in the sense that
$F((\2C\rtimes\2S)_g)\subset(\GLoc_fA)_g$ for all $g\in G$. 
\qed

The following result, which shows that $\Loc_fA^G$ can be computed from $\GLoc_fA$, was
the main motivation for this paper:

\btheor \label{theor-locag}
If $G$ is finite then the functors
\[ \begin{array}{cccc} E: & \Loc_fA^G & \rarr &(\GLoc_fA)^G, \\
   R: & (\GLoc_fA)^G & \rarr & \Loc_fA^G \end{array}\]
are mutually inverse and establish an isomorphism of strict braided tensor categories.
\etheor

\prf By Subsection \ref{ss-ext}, $E: \Loc_fA^G\rarr(\End A)^G$ is a faithful strict tensor
functor, which is full by Proposition \ref{prop-F}. By Proposition \ref{prop-gloc} it
takes its image in $(\GLoc_fA)^G$. By Theorem \ref{theor-E} we have 
$R\circ E=\id_{\Loc_fA^G}$, and $E\circ R=\id_{(\GLoc_fA)^G}$ follows since
$\rho\in(\GLoc_fA)^G$ is the unique right-localized extension to $A_\infty$ of 
$R(\rho)=\rho\restr A^G_\infty$. Therefore $E$ is surjective on objects and thus an
isomorphism. That the braidings of $\Loc_fA^G$ and $(\GLoc_fA)^G$ is clear in view of
their construction. 
\qed

\brem 1. The `size' of $\Loc_fA^G$ will be determined in Corollary \ref{coro-dim}.

2. Clearly the above is a somewhat abstract result, and in concrete models hard work
is required to determine the category $\GLoc_fA$ of twisted representations. (For a
beautiful analysis of orbifolds of affine models in the present axiomatic setting see the
series of papers \cite{x2,lx,klx}.) However, Theorem \ref{theor-locag} can be used to
clarify the structure of $\Loc_fA^G$ quite completely in the holomorphic case, cf.\
Subsection \ref{ss-holom}. 

3. Proposition \ref{prop-gloc} and Theorem \ref{theor-locag} remain true when $G$ is
compact infinite. In order to see this one needs to show that $\2C\rtimes\2S$ is a braided
crossed $G$-category also in the case of infinite $\2S$. In view of the fact that the
existence of $\2C\rtimes\2S$ as rigid tensor category with $G$-action was already
established in \cite{mue06} this can be done by an easy modification of the approach used
in \cite{mue13}. Then the proof of Proposition \ref{prop-F} easily adapts to arbitrary
compact groups.
\erem


\subsection{The equivalence $\Loc_fA^G\rtimes\2S\simeq\GLoc_fA$} \label{ss-equiv2}
Our next aim is to show that the functor $F$ gives rise to an equivalence 
$\Loc_fA^G\rtimes\2S\simeq\GLoc_fA$ of braided crossed $G$-categories. (Even though both 
categories are strict as monoidal categories and as $G$-categories, the functor $F$ will
not be strict.) For the well known definition of a non-strict monoidal functor we refer,
e.g., to \cite{cwm}.

\bprop
If $G$ is finite then the functor $F: \2C\rtimes\2S\rarr\GLoc_fA$ is essentially
surjective, thus a monoidal equivalence.
\eprop

\prf The bulk of the proof coincides with that of \cite[Proposition 3.14]{mue04}, which
remains essentially unchanged. We briefly recall the construction. Pick an interval
$I\in\2K$. Since the $G$-action on $A(I)$ has full spectrum we can find isometries 
$v_g\in A(I), g\in G,$ satisfying 
\[ \sum_g v_gv_g^*=\11, \quad\quad v_g^*v_h=\delta_{g,h}\11, \quad\quad \beta_g(v_h)=v_{gh}. \]
If now $\rho\in\GLoc_fA$ is simple then it is easily verified that
\[ \widetilde{\rho}(\cdot)=\sum_g v_g\,\beta_g\rho\beta_{g^{-1}}(\cdot)\,v_g^* \ \in \GLoc_fA \]
commutes with all $\beta_g$, thus $\widetilde{\rho}\in(\GLoc_fA)^G$. Therefore $\widetilde{\rho}$
restricts to $A^G$, and $\widetilde{\rho}\restr A^G$ is localized in some interval, as was
noted before. In order to show that $\widetilde{\rho}\restr A^G$ is transportable, let $J$ be
some interval, let $\sigma$ be $G$-localized in $J$ and let $s:\rho\rarr\sigma$ be
unitary. Choosing isometries $w_g\in A(J)$ as before and defining $\widetilde{\sigma}$ in
analogy to $\widetilde{\rho}$ and writing $\tilde{s}=\sum_g w_g\beta_g(s)v_g^*$, one easily
verifies that $\tilde{s}$ is a unitary in $\Hom(\widetilde{\rho},\widetilde{\sigma})^G$. Thus
$\widetilde{\rho}\restr A^G_\infty$ is transportable and defines an object of $\Loc_fA^G$. As
in \cite{mue04} one now verifies that $\widetilde{\rho}=E(\widetilde{\rho}\restr A^G)$. Combined
with the obvious fact 
$\rho\prec\widetilde{\rho}$ this implies that every simple object $\rho\in\GLoc_fA$ is a direct
summand of $E(\widetilde{\rho}\restr A^G)=F(\iota(\widetilde{\rho}\restr A^G))$. In view of
Proposition \ref{prop-F} and the fact that $\2C\rtimes\2S$ has splitting idempotents we
conclude that $\rho\cong F(\sigma)$ for some subobject $\sigma$ of
$\iota(\widetilde{\rho}\restr A^G)\in\2C\rtimes\2S$. This implies that $F$ is essentially
surjective, thus an equivalence, which can be made monoidal, see e.g.\ \cite{sr}.
\qed

\brem
In Minkowski spacetimes of dimension $\ge 2+1$, where there are no $g$-twisted
representations, the functor $E$ can be shown to be an equivalence under the weaker
assumption that $G$ is second countable, i.e.\ has countably many irreps, cf.\
\cite{cdr}. Returning to the present one-dimensional situation, it is clear from the
definition of $E$ that $E(\Loc_fA^G\cap\2S')\subset\Loc_fA=(\GLoc_fA)_e$, thus those
$\rho\in\Loc_fA^G$ which satisfy $c_{\rho,\sigma}c_{\sigma,\rho}=\id$ for all
$\sigma\in\2S$ have a localized extension $E(\rho)$. Its seems reasonable to expect that
the restriction of $F$ to the subcategory of $\2C\rtimes\2S$ generated by
$\iota(\2C\cap\2S')$ is an equivalence with $\Loc_fA$ whenever $G$ is second countable. We
have refrained from going into this question this since we are interested in the larger
categories $\Loc_fA^G$ and $\GLoc_fA$, and -- in contradistinction to 
$E: \Loc_fA^G\rarr(\GLoc_fA)^G$ -- the functor $F:\Loc_fA\rtimes\2S\rarr\GLoc_fA$
is almost never essentially surjective (thus an equivalence). The point is that for
$\rho\in\Loc_fA^G$ we have 
$E(\rho)\cong\oplus_i \rho_i$, where the $\rho_i$ are $g_i$-localized and the $g_i$
exhaust a whole conjugacy class since $E(\rho)$ is $G$-invariant. Since the direct sum is
finite, we see that the image of $E:\2C\rtimes\2S\rarr \GLoc_fA$ can contain only objects
$\sigma$ whose degree $\del\sigma$ belongs to a finite conjugacy class. Since `most' infinite
non-abelian compact groups have infinite conjugacy classes, $F$ will in general not be
essentially surjective. (At least morally this is related to the fact \cite{km} that the
quantum double of a compact group $G$ admits infinite dimensional irreducible
representations whenever $G$ has infinite conjugacy classes.) If, on the other hand, we
consider $E(\rho)$ where $d(\rho)=\infty$, the analysis of $E(\rho)$ becomes considerably
more complicated.  
\erem

\bcoro \label{coro-dim}
Under the assumptions of Theorem \ref{theor-M2} we have
\[ \dim\Loc_fA^G=|G|\,\dim\GLoc_fA. \]
\ecoro

\prf Follows from $\GLoc_fA\cong\Loc_fA^G\rtimes\2S$ and 
$\dim\2C\rtimes\2S=\dim\2C/\dim\2S=\dim\2C/|G|$, cf.\ \cite{mue06}.
\qed

In order to prove the equivalence $\GLoc_fA\simeq\Loc_fA^G\rtimes\2S$ of braided crossed
$G$-categories we need to consider the $G$-actions and the braidings. For the general
definition of functors of $G$-categories we refer to \cite{tam}, see also \cite{carr} and
the references given there. Since our categories are strict as tensor categories and as
$G$-categories, i.e.\  
\begin{equation}
\ba{rcll} \gamma_{gh}(X) &=& \gamma_g\circ\gamma_h(X) & \forall g,h,X, \\
  \gamma_g(X\otimes Y) &=& \gamma_g(X)\otimes\gamma_g(Y) & \forall g,X,Y,\ea
\label{tur}\end{equation}
we can simplify the definition accordingly:

\newarrow{Congruent} 33333
\bdefin \label{def-Gcat}
A functor $F: \2C\rarr\2C'$ of categories with strict actions $\gamma_g,\gamma'_g$ of a
group $G$ is a functor together with a family of natural isomorphisms 
$\eta(g):F\circ\gamma_g\rarr\gamma'_g\circ F$ such that 
\[ \begin{diagram} 
  F\circ\gamma_{gh}(X) && \rTo^{\eta(gh)_X} && \gamma'_{gh}\circ F(X) \\
  \dCongruent &&&& \dCongruent \\
  F\circ\gamma_g\circ\gamma_h(X) & \rTo_{\eta(g)_{\gamma_h(X)}} & 
   \gamma'_g\circ F\circ\gamma_h(X) & \rTo_{\gamma'_g(\eta(h)_X)} & 
  \gamma'_g\circ\gamma'_h\circ F(X) \end{diagram} \]
commutes. (There is no further condition on $F$ if $\2C,\2C',\gamma,\gamma'$ are monoidal.)

A functor of braided crossed $G$-categories is a monoidal functor of $G$-categories that
respects the gradings and satisfies $F(c_{X,Y})=c_{F(X),F(Y)}$ for all $X,Y\in\2C$. 
\edefin

\btheor \label{theor-M2}
Let $A=(\2H_0,A(\cdot),\Omega)$ be as before and $G$ finite. Then  
\[ F: \ \2C\rtimes\2S\ \rarr \ \GLoc_fA \]
is an equivalence of braided crossed $G$-categories.
\etheor

\prf It only remains to show that $F$ is a functor of $G$-categories and that it preserves
the braidings. Let $(\rho,p)\in\2C\rtimes\2S$. Then $\beta_g((\rho,p))=(\rho,\beta_g(p))$, 
where $\beta_g(p)$ is the obvious $G$-action on $\2C\rtimes_0\2S$. Recall that
$F((\rho,p))\in\End\,A_\infty$ was defined as $v_{(\rho,p)}^*E(\rho)(\cdot)v_{(\rho,p)}$,
where $v_{(\rho,p)}\in A_\infty$ satisfies $v_{(\rho,p)}v^*_{(\rho,p)}=E(p)$. (For $p=1$
we choose $v_{(\rho,p)}=1$.) Since $E(\rho)$ commutes with $\gamma_g$ we have 
$\gamma_g(F((\rho,p)))=\gamma_g(v_{(\rho,p)})^*E(\rho)(\cdot)\gamma_g(v_{(\rho,p)})$.
Because of $\gamma_g(v_{(\rho,p)})\gamma_g(v_{(\rho,p)})^*=\gamma_g(p)$, the isometries
$\gamma_g(v_{(\rho,p)})$ and $v_{\beta_g(\rho,p)}$ have the same range projection. Thus 
$\eta(g)_{(\rho,p)}=\gamma_g(v_{(\rho,p)})^*v_{(\rho,\beta_g(p))}$ is unitary and one
easily verifies
$\eta(g)_{(\rho,p)}\in\Hom(F\circ\beta_g(\rho,p),\gamma_g\circ F(\rho,p))$ as well as the
commutativity of the above diagram.

It remains to show that the functor $F$ preserves the braidings. We first show that
$F(c_{\rho,\sigma})=c_{F(\rho),F(\sigma)}$ holds if $\rho,\sigma\in\2C=\Loc_fA^G$. By
Theorem \ref{theor-E}, $E(\rho), E(\sigma)$ are $G$-invariant, thus by the $G$-covariance
of the braiding we have $c_{E(\rho),E(\sigma)}\in A^G_\infty$. Thus the braiding of
$E(\rho),E(\sigma)$ as constructed in Section \ref{s-crossed} restricts to a braiding of
$\rho,\sigma$ and by uniqueness of the latter this restriction coincides with
$c_{\rho,\sigma}$. Thus $c_{E(\rho),E(\sigma)}=E(c_{\rho,\sigma})$ as claimed.
The general result now is an obvious consequence of the naturality of the braidings of
$\2C\rtimes\2S$ and of $\GLoc_fA$ together with the fact that every object of
$\2C\rtimes\2S$ and of $\GLoc_fA$ is a subobject of one in $\2C$ and $(\GLoc_fA)^G$,
respectively. 
\qed


\section{Orbifolds of completely rational chiral CFTs} \label{s-cplrtl}
\subsection{General theory}
So far, we have considered an arbitrary QFT $A$ on $\7R$ subject to the technical
condition that also $A^G$ be a QFT on $\7R$, some of the results assuming finiteness of 
$G$. The situation that we are really interested is the one where $A$ derives from a
chiral QFT on $S^1$ by restriction to $\7R$. Recall that in that case $\Loc_{(f)}A^G$ has
a `physical' interpretation as a category $\Rep_{(f)}A$ of representations.

\bprop Let $A$ be a completely rational chiral QFT with finite symmetry group $G$. Then
the restrictions to $\7R$ of $A$ and $A^G$ are QFTs on $\7R$.
\eprop

\prf In view of the discussion in Subsection \ref{ss-chiral} it suffices to know that the
chiral orbifold theory $A^G$ on $S^1$ satisfies strong additivity. In \cite{x2} it was
proven that finite orbifolds of completely rational chiral QFTs are again completely
rational, in particular strongly additive.
\qed

Applying the results of \cite{klm} we obtain:

\btheor \label{theor-M3}
Let $(\2H_0,A,\Omega)$ be a completely rational chiral CFT and $G$ a finite
symmetry group. Then the braided crossed $G$-category $\GLoc_fA$ has full $G$-spectrum,
i.e.\ for every $g\in G$ there is an object $\rho\in\GLoc_fA$ such that
$\del\rho=g$. Furthermore, for every $g\in G$ we have
\[ \sum_{\rho\in(\GLoc_fA)_g} (\dim\rho)^2 =\sum_{\rho\in\Rep_fA} (\dim\rho)^2=\mu(A), \]
where the sums are over the the equivalence classes of irreducible objects of degree $g$
and $e$, respectively.
\etheor

\prf
By \cite{x2}, the fixpoint theory $A^G$ is completely rational, thus by \cite{klm} the 
categories $\Rep_fA^G\cong\Loc_fA^G$ are modular. Now,
$\GLoc_fA\cong\Loc_fA^G\rtimes\2S$, and fullness of the $G$-spectrum follows by
\cite[Corollary 3.27]{mue13}. The statement on the dimensions follows from
\cite[Proposition 3.23]{mue13}.
\qed

\brem 1. It would be very desirable to give a direct proof of the fullness of the $G$-spectrum
of $\GLoc_fA$ avoiding reference to the orbifold theory $A^G$ via the equivalence
$\GLoc_fA\simeq \Loc_fA^G\rtimes\2S$. This would amount to showing directly that
$g$-localized transportable endomorphisms of $A_\infty$ exist for every $g\in G$. Since
our proof relies on the fairly non-trivial modularity result for $\Loc_fA^G$, cf.\
\cite{klm} together with \cite{x2}, this might turn out difficult. 

2. In the VOA setting, Dong and Yamskulna \cite{dy} have shown that there exist twisted 
representations for all $g\in G$. Since \cite[Proposition 3.23]{mue13} is a purely
categorical result, the above conclusion also holds in the VOA setting as soon as one can
establish that the G-twisted representations form a rigid tensor category.

3. It may be useful to summarize the situation in a diagram:
\begin{diagram}
 \Loc_fA & \subset & \ \ \GLoc_fA  \\
 \dTo^{\2C\leadsto\2C^G}\uTo & & \dTo\uTo_{\2C\leadsto\2C\rtimes\2S} \\
 \Loc_fA^G\cap\2S'\ \  & \subset & \Loc_fA^G
\end{diagram}
The horizontal inclusions are full, $\Loc_fA$ being the degree zero subcategory of
$\GLoc_fA$. If $G$ is abelian, the $G$-grading passes to $\Loc_fA^G$ (see \cite{mue13})
and $\Loc_fA^G\cap\2S'$ is its degree zero subcategory. Moving from left to right or from
top to bottom, the dimension of the categories are multiplied by $|G|$. In the upper line
this is due to Theorem \ref{theor-M3} and in the lower due to the results of
\cite{mue11}. Together with $\dim\2C=|G|\cdot\dim\2C\rtimes\2S$ this implies
$\dim\Loc_fA^G=|G|^2\dim\Loc_fA$, as required by \cite{klm}. (In fact, this latter
identity together with \cite[Proposition 3.23]{mue13} provides an alternative proof of the
completeness of the $G$-spectrum of $\GLoc_fA$.) Furthermore, the upper left and lower
right categories are modular, whereas $\Loc_fA^G\cap\2S'$ is not (whenever
$G\ne\{e\}$). The passage $\Loc_fA^G\cap\2S' \ \leadsto \ \Loc_fA$ is the `modular
closure' from \cite{mue06,brug} and $\Loc_fA^G\cap\2S' \ \leadsto \ \Loc_fA^G$ is the
`minimal modularization', conjectured to exist for every premodular category, cf.\ \cite{mue11}.
\erem

We briefly discuss the modularity of $\GLoc_fA$. In \cite{t2}, a braided crossed
$G$-category $\2C$ was called modular if its braided degree zero subcategory $\2C_e$ is
modular in the usual sense \cite{t1}. This definition seems somewhat unsatisfactory since
it does not take the nontrivially graded part of $\2C$ into account. In \cite{kir3}, the
vector space 
\[ \2V_\2C=\bigoplus_{i\in I}\bigoplus_{g\in G} \Hom(\beta_g(X_i),X_i), \]
where $I$ indexes the isomorphism classes of simple objects in $\2C$, is introduced and an
endomorphism $S\in\End\,\2V_\2C$ is defined by its matrix elements
\[ S((X,u),(Y,v)))=\quad \begin{tangle}
\hcoev\step\hcoev \\ 
\id\step\O{u}\step\id\step\id\\
\id\step\hxx\step\id\\
\id\step\O{v}\step\id\step\id\\
\id\hstep\mobj{X}\hstep\hxx\mobj{Y}\step\id\\
\hev\step\hev
\end{tangle}\]
where $\del X=g,\ \del Y=h$ and $u: \beta_h(X)\rarr X,\ v: \beta_g(Y)\rarr Y$. A braided
$G$-crossed fusion category is modular (in the sense of \cite{kir3}) if the endomorphism
$S$ is invertible. 

\bprop 
Let $(\2H_0,A,\Omega)$ be a completely rational chiral CFT and $G$ a finite
symmetry group. Then the braided crossed $G$-category $\GLoc_fA$ is modular in the sense
of \cite{kir3}.
\eprop

\prf As used above, the braided categories $\Loc_fA=(\GLoc_fA)_e$ and 
$\Loc_fA^G\cong(\GLoc_fA)^G$ are modular. Now the claim follows by \cite[Theorem
10.5]{kir3}. 
\qed

The preceding discussions have been of a very general character. In the next subsection
they will be used to elucidate completely the case of holomorphic orbifolds, where our
results go considerably beyond (and partially diverge from) those of \cite{dvvv}. In the 
non-holomorphic case it is clear that comparably complete results cannot be hoped for. 
Nevertheless already a preliminary analysis leads to some surprising results and
counterexamples, cf.\ the final subsection.


\subsection{Orbifolds of holomorphic models} \label{ss-holom}
\bdefin
A holomorphic chiral CFT is a completely rational chiral CFT with trivial representation 
category $\Loc_fA$. (I.e.,\ $\Loc_fA$ is equivalent to $\mathrm{Vect}_f\7C$.)
\edefin 

\brem By the results of \cite{klm}, a completely rational chiral CFT is holomorphic iff
$\mu(A)=1$ iff $A(E')=A(E)'$ whenever $E=\cup_{i=1}^n I_i$ where $I_i\in\2I$ with mutually
disjoint closures.
\erem

\bcoro \label{coro-holo}
Let $A$ be a holomorphic chiral CFT acted upon by a finite group $G$. Then $\GLoc_fA$ has 
precisely one isomorphism class of simple objects for every $g\in G$, all of these objects
having dimension one.
\ecoro

\prf By Theorem \ref{theor-M3}, we have $\dim(\GLoc_fA)_g=1$ for all $g\in G$. Since the
dimensions of all objects are $\ge 1$, the result is obvious.
\qed

\brem 1. In \cite{mue05}, where the invertible objects of $\GLoc_fA$ were called soliton 
automorphisms, it is shown that these objects can be studied in a purely local manner.

2. Let $A$ be a holomorphic chiral CFT, and pick an interval $I\in\2K$. By Corollary
\ref{coro-holo} there is just one isoclass of simple objects in $(\GLoc_fA)_g$ for every
$g\in G$. Since the objects of $\GLoc_fA$ are transportable endomorphisms of $A_\infty$,
we can pick, for every $g\in G$, representer $\rho_g$ that is $g$-localized in $I$. By
Lemma \ref{l-br1}, $\rho_g$ restricts to an automorphism of $A(I)$. Furthermore, we can
choose unitaries $u_{g,h}\in\Hom_{A(I)}(\rho_g\rho_h,\rho_{gh})$. In other words,
we have a homomorphism 
\[ G\rarr\Aut A(I)/\mathrm{Inn} A(I)=:\mathrm{Out}A(I), \quad g\mapsto[\rho_g], \]
thus a `G-kernel', cf.\ \cite{suth}. We recall some well known facts: The associativity 
$(\rho_g\rho_h)\rho_k=\rho_g(\rho_h\rho_k)$ implies the existence of
$\alpha_{g,h,k}\in\7T$ such that 
\[ u_{gh,k}u_{g,h}=\alpha_{g,h,k}\,u_{g,hk}\,\rho_g(u_{h,k}) \quad\forall g,h,k. \]
A tedious but straightforward computation using four $\rho$'s shows that 
$\alpha:G\times G\times G\rarr\7T$ is a 3-cocycle, whose cohomology class 
$[\alpha]\in H^3(G,\7T)$ does not depend on the choice of the $\rho$'s and of the
$u$'s. Thus $[\alpha]$ is an obstruction to the existence of representers $\rho_g$ for
which $g\mapsto\rho_g$ is a homomorphism $G\rarr\Aut A(I)$. (Actually, since in QFT the
algebras $A(I)$ are type III factors with separable predual, the converse is also true: If
$[\alpha]=0$ then one can find a homomorphism $g\mapsto\rho_g$, cf.\ \cite{suth}.) 
\erem

\bnix
For a further analysis it is more convenient to adopt a purely categorical viewpoint. 
Starting with the category $\GLoc_fA$ of a holomorphic theory $A$, we don't lose any
information by throwing away the non-simple objects and the zero morphisms. In this way we
obtain a categorical group $\2C$, i.e.\ a monoidal groupoid where all objects have a
monoidal inverse. The set of isoclasses is the group $G$. In the general $k$-linear case
it is well known that such categories are classified up to equivalence by
$H^3(G,k^*)$. This is shown by picking an equivalent skeletal tensor category
$\widetilde{\2C}$, i.e.\ a full subcategory with one object per isomorphism class. Even if
$\2C$ is strict, $\widetilde{\2C}$ in general is not, and the associativity constraint
defines an element of $H^3(G,k^*)$. It is thus clear that 3-cocycles on $G$ will play a
r\^ole in the classification of the braided crossed $G$-categories associated with
holomorphic QFTs. In view of \cite{dvvv,dpr,dw} and \cite{dm,dy} this is hardly
surprising. Yet, the situation is somewhat more involved than anticipated by most authors
since a classification of the possible categories $\GLoc_fA$ -- and therefore of the
categories $\Loc_fA^G$ -- must also take the $G$-action on $\GLoc_fA$ and the braiding
into account. 

If one considers {\it braided} categorical groups, $G$ must be abelian and one has a 
classification in terms of $H^3_{\mathrm{ab}}(G,k^*)$, cf.\
\cite{js}. ($H^3_{\mathrm{ab}}(G,k^*)$ is Mac Lane's cohomology \cite{ml} for abelian
groups.) The requirement that $G$ be abelian disappears if one admits a non-trivial
$G$-action and considers braided crossed $G$-categories. One finds \cite{t2} that
(non-strict) skeletal braided crossed $G$-categories with strict $G$-action in the sense
of (\ref{tur}) are classified in terms of Ospel's quasiabelian cohomology
$H^3_{\mathrm{qa}}(G,k^*)$ \cite{osp}. Unfortunately, this is still not sufficient for our
purposes. Namely, assume we have a braided crossed $G$-category $\2C$ that is also a
categorical group (and thus a categorical $G$-crossed module in the sense of
\cite{carr}). Even if $\2C$ is strict monoidal and satisfies (\ref{tur}) -- as 
our categories $\GLoc\,A$ and $\2C\rtimes\2S$ do -- an equivalent skeletal category
$\widetilde{\2C}$ in general will not satisfy (\ref{tur}). It is clear that for a completely
general classification of braided crossed $G$-categories that are categorical groups one
can proceed along similar lines as in the classifications cited above. We will supply the
details in the near future \cite{mue17}, also elucidating the r\^ole of the twisted
quantum doubles $D^\omega(G)$ \cite{dpr} in the present context. (Note that the modular
category $D^\omega-\Mod$ contains the symmetric category $G-\Mod$ as a full subcategory,
and $D^\omega-\Mod\rtimes G-\Mod$ is a braided crossed $G$-category with precisely one 
invertible object of every degree. However, not every such category is equivalent to 
$D^\omega-\Mod\rtimes G-\Mod$ for some $[\omega]\in H^3(G,\7T)$!)
\enix


\subsection{Some observations on non-holomorphic orbifolds}
In the previous subsection we have seen that a holomorphic chiral CFT $A$ has (up to
isomorphism) exactly one simple object of degree $g\in G$, and this object has dimension
one, thus is invertible. This allows a complete classification of the categories $\GLoc A$
and $\Loc A^G\cong(\GLoc A)^G$ that can arise.

It is clear that in the non-holomorphic case ($\Loc_fA\not\simeq\mathrm{Vect}_\7C$) there
is no hope of obtaining results of this completeness. The best one could hope for would be
a classification of the categories $\GLoc_fA$ that can arise from CFTs with prescribed
$\Loc_fA\simeq(\GLoc_fA)_e$, but for the time being this is far out of reach. We therefore
content ourselves with some comments on a more modest question. To wit, we ask whether a 
non-holomorphic completely rational CFT $A$ admits {\it invertible} $g$-twisted
representations  for every $g\in G$. (As we have seen, this is the case for holomorphic
$A$.) It turns out that the existence of a braiding (in the sense of crossed
$G$-categories) provides an obstruction:  

\blemma
Let $\2C$ be a braided crossed $G$-category. If there exists an invertible object of
degree $g\in G$ then 
\[ \gamma_g(X)\cong X \quad \forall X\in\2C_e. \]
\elemma

\prf Let $X\in\2C_e$ and $Y\in\2C_g$. Then the braiding gives rise to isomorphisms
$c_{X,Y}: X\otimes Y\rarr Y\otimes X$ and $c_{Y,X}: Y\otimes X\rarr\gamma_g(X)\otimes Y$.
Composing these we obtain an isomorphism $X\otimes Y \rarr \gamma_g(X)\otimes Y$. If $Y$
is invertible, we can cancel it by tensoring with $\ol{Y}$, obtaining the desired
isomorphism $X\stackrel{\cong}{\longrightarrow}\gamma_g(X)$. 
\qed

\bcoro \label{coro1}
Let $\2C$ be a braided crossed $G$-category and let $g\in G$. If there exists $X\in\2C_e$
such that $\gamma_g(X)\not\cong X$ then there exists no invertible $Y\in\2C_g$. 
\ecoro

\brem The condition $\gamma_g(X)\cong X\ \forall X\in\2C_e$ is necessary in order for
the existence of invertible objects of degree $g$, but of course not sufficient. In any
case, there are many chiral CFTs where the corollary, as applied to $\GLoc_fA$, excludes
invertible $g$-twisted representations for $g\ne e$. One such class will be considered below.
\erem

We apply the above results to the $n$-fold direct product $A=B^{\otimes n}$ of a
completely rational chiral CFT $B$, on which the symmetric group $S_n$ acts in the obvious
fashion. We first note that every irreducible $\pi\in\Rep_fA$ is unitarily
equivalent to a direct product $\pi_1\otimes\cdots\otimes\pi_n$ of irreducible
$\pi_i\in\Rep_fB$, cf.\ \cite{klm}. Thus the equivalence classes of simple objects of
$\Loc_fA$ are the $n$-tuples of equivalence classes of simple objects of $\Loc_fB$, and
$S_n$ acts on them by permutation.  

\bcoro \label{coro2}
Let $B$ be a completely rational chiral CFT and let $n\ge 2$. Consider $A=B^{\otimes n}$
with the permutation action of $G=S_n$. If $B$ is not holomorphic then $\GLoc_fA$ contains
no invertible object $\rho$ with $\del\rho\ne e$. 
\ecoro

\prf Since $B$ is not holomorphic we can find a simple object $\sigma\in\Loc_fB$ such that
$\sigma\not\cong\11$. If $g\in S_n$ with $g\ne e$ there is $i\in\{1,\ldots,n\}$ such that
$g(i)\ne i$. Consider an object $\rho=(\rho_1,\ldots,\rho_n)\in\Loc_fA$ where $\rho_i=\11$
and $\rho_{g(i)}=\sigma$. Now it is clear that $\gamma_g(\rho)\not\cong\rho$, and
Corollary \ref{coro1} applies.
\qed

For any tensor category $\2C$ we denote by $\Pic(\2C)$ the full monoidal subcategory of
invertible objects. (In a $*$-category these are precisely the objects of dimension one.) 

\bcoro \label{coro3}
Let $B$ be a completely rational chiral CFT. Consider $A=B^{\otimes n}$ for $n\ge 2$ and
let $G\subset S_n$ be a subgroup. If $B$ is non-holomorphic then
\[ \Pic(\Loc_fA^G)\cong\Pic\left((\Loc_fA)^G\right). \]
\ecoro

\prf We may assume $G\ne\{e\}$ since otherwise there is nothing to prove. By Theorem
\ref{theor-locag} we have $\Loc_fA^G\cong(\GLoc_fA)^G$. Let now
$\rho\in\Pic(\Loc_fA^G)$. Then $E(\rho)\in\Pic((\GLoc_fA)^G)$, and by 
Corollary \ref{coro2} we have $\del E(\rho)=e$, thus $E(\rho)\in\Pic((\Loc_fA)^G)$.
The rest follows as in Subsection \ref{ss-equiv1}. 
\qed

Thus, in permutation orbifold models, the Picard category $\Pic(\Loc_fA^G)$ is
determined already by $\Pic(\Loc_fA)$ and the $G$-action on it, i.e.\ we do not need to
know the $g$-twisted representations of $A$ for $g\ne e$. We recall that a subgroup
$G\subset S_n$ is called transitive if for each $i,j\in\{1,\ldots,n\}$ there exists
$g\in G$ such that $g(i)=j$.

\bcoro \label{coro4}
Let $B$ be a non-holomorphic completely rational chiral CFT. Consider $A=B^{\otimes n}$
for $n\ge 2$ and let $G\subset S_n$ be a transitive subgroup. Then the isomorphism classes
in $\Pic(\Loc_fA^G)$ are in 1-1 correspondence with the pairs $([\sigma],\lambda)$, where
$[\sigma]$ is an isomorphism class in $\Pic(\Loc_f B)$ and
$\lambda\in\widehat{G}_1=\widehat{G_{ab}}$ is a one-dimensional character of $G$. 
\ecoro

\prf Let $\rho$ be an invertible object of $\Loc_fA^G$. By Corollary \ref{coro3}, we
have $E(\rho)\cong(\sigma_1,\ldots,\sigma_n)$ where the $\sigma_i$ are invertible objects
of $\Loc_fB$. By Subsection \ref{ss-ext}, $E(\rho)$ is invariant under the $G$-action on
$\Loc_fA$, and since the latter transitively permutes the $\sigma_i$ there is
$\sigma\in\Pic(\Loc_fB)$ such that $\sigma_i\cong\sigma$ for all $i$. Now, by \ref{nix-S}
we know that for every $\lambda\in\widehat{G}_1$ there exist localized unitaries
$u_\lambda\in A_\infty$ such that $\beta_g(u_\lambda)=\lambda(g)u_\lambda$. In restriction
to $A^G_\infty$, the localized isomorphisms $\mathrm{Ad}\,u_\lambda$ are inequivalent
invertible objects $\rho_\lambda\in\Pic(\Loc_fA^G)$. Now the claimed bijection follows by
picking one representer $\sigma$ for each isoclass $[\sigma]$ in $\Pic(\Loc_fB)$ and
mapping $([\sigma],\lambda)\mapsto[(\sigma,\ldots,\sigma)\otimes\rho_\lambda]$.
\qed

\brem At this place in the preceding version of this paper, which will appear in
Commun. Math. Phys., I claimed that the results of this subsection are in contradiction
to what can be derived from certain statements in \cite{ban2}, which in turn follow from
\cite{ban1}. This claim was wrong, being based on an erroneous deduction from the
statements in \cite{ban1,ban2}. I regret this mistake. In fact, Bantay has provided me
with a convincing argument to the effect that also his completely independent methods
imply Corollary \ref{coro4} above. His argument relies on the formula 
\cite[eq.\ (15)]{ban1} for the $S$-matrix of the permutation orbifold, which can be 
traced back to the character formula \cite[eq.\ (5)]{ban1}. 

However, I remain unconvinced by the justification of the latter given in \cite{ban1} and
still recommend \cite{ml2}, where a vigorous case is made for rigorous proof in
theoretical physics. (As to the labelling of the irreducible sectors of the permutation
orbifold stated  in \cite{ban1} without even a hint of proof, such a proof has recently
been provided in \cite{klx}.) 
\erem

\noindent{\it Acknowledgments.} The research reported here was presented at the workshop
`Tensor Categories in Mathematics and Physics' which took place at the Erwin Schr\"odinger
Institute, Vienna, in June 2004. I am grateful to the ESI for hospitality and financial 
support and to the organizers for the invitation to a very stimulating meeting.


\end{document}